\documentclass[a4paper,12pt]{article}

\usepackage[]{setspace}
\usepackage{subcaption}
\usepackage{lineno,hyperref}
\usepackage{fullpage}
\usepackage{comment}
\usepackage{xcolor}
\usepackage{amsmath,amssymb,amsfonts}
\usepackage{graphicx}
\usepackage{tikz}
\usepackage{pgfplots}
\usepackage{algorithm}
\usepackage{dsfont}
\usepackage{algpseudocode}
\usetikzlibrary{calc}
\usepackage{siunitx}
\usepackage{listings}
\usepackage{pgf-umlcd}
\usepackage{pgf-umlsd}
	\colorlet{punct}{red!60!black}
\definecolor{background}{HTML}{EEEEEE}
\definecolor{delim}{RGB}{20,105,176}
\colorlet{numb}{magenta!60!black}
\lstdefinelanguage{json}{
	basicstyle=\normalfont\ttfamily,
	numbers=left,
	numberstyle=\scriptsize,
	stepnumber=1,
	numbersep=8pt,
	showstringspaces=false,
	breaklines=true,
	frame=lines,
	backgroundcolor=\color{background},
	literate=
	*{0}{{{\color{numb}0}}}{1}
	{1}{{{\color{numb}1}}}{1}
	{2}{{{\color{numb}2}}}{1}
	{3}{{{\color{numb}3}}}{1}
	{4}{{{\color{numb}4}}}{1}
	{5}{{{\color{numb}5}}}{1}
	{6}{{{\color{numb}6}}}{1}
	{7}{{{\color{numb}7}}}{1}
	{8}{{{\color{numb}8}}}{1}
	{9}{{{\color{numb}9}}}{1}
	{:}{{{\color{punct}{:}}}}{1}
	{,}{{{\color{punct}{,}}}}{1}
	{\{}{{{\color{delim}{\{}}}}{1}
	{\}}{{{\color{delim}{\}}}}}{1}
	{[}{{{\color{delim}{[}}}}{1}
	{]}{{{\color{delim}{]}}}}{1}
	{``}{\textquotedblleft}1,
}
\modulolinenumbers[5]

\newcommand{\Fluid}{\mathcal{F}} 
\newcommand{\Real}{\mathbb{R}} 
\newcommand{\Alemap}{\mathcal{A}} 
\newcommand{\ALE}{ALE} 
\newcommand{\Vel}{u} 
\newcommand{\Pres}{p} 
\newcommand{\Density}{\rho} 
\newcommand{\Viscosity}{\mu} 
\newcommand{\tvel}{U} 
\newcommand{\angvel}{\omega} 
\newcommand{\Rmat}{R} 
\newcommand{\Angle}{\theta}
\newcommand{\Inertia}{I} 
\newcommand{\mass}{m} 
\newcommand{\CenterMass}{x^{CM}} 
\newcommand{\Solid}{\mathcal{S}} 
\newcommand{\normal}{n} 

\newlength{\dhatheight}

\begin{document}

\begin{center}
    \large
    \textbf{Mathematical and computation framework for moving and colliding rigid bodies in a Newtonian fluid}
        
    \vspace{1.0cm}
	\normalsize
    \textbf{Céline Van Landeghem, Luca Berti, Vincent Chabannes, \\ Christophe Prud'homme}\\
	\textit{Cemosis, IRMA UMR 7501, CNRS, Université de Strasbourg, France}

	\vspace{0.3cm}

	\textbf{Laetitia Giraldi}\\
	\textit{CALISTO team, INRIA, Université Côte d'Azur, France}

	\vspace{0.3cm}

	\textbf{Agathe Chouippe, Yannick Hoarau}\\
	\textit{ICube Laboratory, UMR 7357, Université de Strasbourg, France}

    \vspace{1.0cm}
    \textbf{Abstract}
\end{center}

\noindent We studied numerically the dynamics of colliding rigid bodies in a Newtonian fluid.	
The finite element method is used to solve the fluid-body interaction and the fluid motion 
is described in the Arbitrary-Lagrangian-Eulerian framework. To model the interactions between bodies, 
we consider a repulsive collision-avoidance model, 
defined by R. Glowinski in \cite{glowinski_linite_nodate}.
The main emphasis in this work is the generalization of this collision model 
to multiple rigid bodies of arbitrary shape. Our model first uses a narrow-band 
fast marching method to detect the set of colliding bodies.  
Then, collision forces and torques 
are computed for these bodies via a general expression, which does not depend on their shape. 
Numerical experiments examining the performance of the narrow-band fast marching method 
and the parallel execution of the collision algorithm are discussed. 
We validate our model with literature results and show various applications 
of colliding bodies in two and three dimensions. In these applications, 
the bodies either move due to gravity, a flow, or can actuate themselves. Finally, we 
present a tool to create arbitrary shaped bodies in complex already discretized 
fluid domains, enabling conforming body-fluid interface and allowing to perform 
simulations of fluid-body interactions with collision treatment in these realistic 
environments. 
All simulations are conducted with the \textbf{Feel++} open source library.

\vspace{0.3cm}


%
%
\section*{Introduction}

Fluid flows laden with particles are common in 
industrial and biological processes, such as fluidization, the cell transport in arteries, 
or the simulated motion of articulated micro-swimmers in the human body. 
Due to the particle volume and the confined environments, these processes are characterized by inter-particle 
interactions. The modeling of these interactions, based on collision detection algorithms, 
and the computation of lubrication and collision forces, is challenging for arbitrary shaped particles  
and add further complexity to the coupled fluid-solid 
interaction problem.\newline 

In recent years, various approaches simulating solid-solid or solid-wall interactions 
have been developed and proposed in the literature. These approaches are based on 
collision and lubrication forces, added to the system when solving the fluid-solid 
interactions \cite{wachs_modeling_2023}.
When the distance between solid surfaces is small, approximately the mesh size of 
the computational grid, the fluid placed between the solids is squeezed out and the 
hydrodynamics are not resolved. The resulting underestimation of the hydrodynamic forces 
is compensated by a short-range repulsive correction given analytically and 
called lubrication force. In \cite{jain_collision_2019}, this lubrication force is defined by a 
constant expression and in \cite{singh_distributed_2003} its magnitude is inversely proportional to the 
distance between the solid surfaces. In some approaches, the lubrication model 
allows avoiding the direct contact between the solid surfaces \cite{glowinski_linite_nodate}. 
Such contact avoidance schemes are frequently used in literature to 
simulate the interactions between spherical bodies. In this paper, we will present 
the generalization of these schemes for complex shaped and articulated bodies.  
The contact avoidance models often require the adjustment of stiffness parameters, for which 
the optimal values are not known. Over- or 
underestimation of these values can falsify the results. To overcome 
this issue, the authors in \cite{broms_barrier_2023} solve a minimization problem to determine the 
minimal force magnitude ensuring a fixed separation distance between 
the solids at each time step and in \cite{lefebvre_numerical_2009} the authors introduce the gluey particle method.\newline 

A collision force is applied, when the lubrication correction does not avoid 
the direct contact between the solids. Common collision forces are based on the 
soft-sphere or hard-sphere approach. 
When two solids are in direct contact or overlapping, the soft-sphere approach 
defines the collision force according to the geometrical characteristics of this 
contact \cite{ardekani_numerical_2016}. The approach requires 
a small time step to resolve the collision process.    
In contrast, the hard-sphere approach \cite{jain_collision_2019} assumes that the collision force is impulsive. 
During the collision process, the velocities of the interacting solids change 
instantaneously at the time of contact. 
The magnitude of the new velocities depends on the pre-collision velocities 
and some physical parameters.  \newline 

For large fluid systems, efficient collision detection algorithms are mandatory. 
These algorithms identify the pairs of solids that are actually interacting, 
before computing the collision forces only for the respective pairs, which reduces 
the computational costs. 
Most collision detection algorithms contain two phases, the broad-phase and the narrow-phase.
First, the broad-phase uses spatial partitioning or sorting methods to determine 
the smallest possible set of neighbouring pairs that are likely to interact. 
In spatial partitioning, the fluid system is divided into regions and all solids 
assigned to the same region are considered as neighbours. The data structures employed 
for spatial partitioning and neighbour detection are lists \cite{muth_collision_nodate}, trees \cite{vemuri_efficient_1998}, adapted for 
large number of solids, or hash tables \cite{pouchol_hierarchical_nodate}, efficient for solids with a wide range of sizes. 
The spatial sorting methods, as the Sort and Sweep algorithm \cite{real_time}, sort the solids in space 
to determine overlapping. Then, the narrow-phase detects whether these pairs are 
actually in contact by computing the distance between the solid surfaces. 
The distance computation being evident for spherical solids, it is a 
challenging and often expensive task for complex shaped or articulated bodies. 
In the literature, multiple analytical methods for specific geometries or 
general algorithms for arbitrary solid shapes are proposed \cite{wachs_modeling_2023}.
Our model uses the fast marching method \cite{sethian_fast_1996} to track the distance between any pair 
of arbitrary shaped bodies. Focusing on fluid systems containing up to one hundred solids, the 
implementation of a broad-phase algorithm is not necessary. However, to avoid the high 
computational costs associated to the fast marching method, we adapt a narrow-band 
approach, executable in parallel.  \newline 

The paper is organized as follows. 
Section \ref{rigidmotion} details the coupled problem describing the 
fluid-solid interaction.  The collision model including the collision detection 
algorithm based on the narrow-band fast marching method and the computation of the 
contact avoidance lubrication force is presented in section \ref{collisionmodel}.
The performance verification and the validation of the model are done in sections \ref{experiments} and \ref{validation}.
Section \ref{applications} shows some applications of the collision model, 
including simulations of single complex shaped and articulated solids as well as 
spherical multi-bodies interactions. A summary and conclusions are given in section 
\ref{conclusion}. \newline 

\section{Rigid body moving in a fluid} \label{rigidmotion}

This section briefly describes the coupled fluid-rigid body interaction 
problem. A detailed description is given in \cite{berti_fluid-rigid_nodate}.

\paragraph*{The fluid model} 
The fluid model considered throughout the paper is incompressible and Newtonian, 
and the motion of the fluid domain is modelled using the Arbitrary-Lagrangian-Eulerian 
(\ALE) formalism \cite{formaggia}. Let $\Fluid_t \subseteq \Real^d, 
d\in\{2,3\},$ denote the domain occupied by the 
fluid at time $t$, where $t\in [0,T]$ and $T$ the final time of the simulation. 
Let $\Alemap_t: \Fluid_0 \to \Fluid_t$ be the 
\ALE\ map mapping the reference fluid domain $\Fluid_0$ to the current domain $\Fluid_t$, 
defined as $\Alemap_t(X)=\Alemap(t,X)=X + x(t,X)$ with $x(t,X)$ the displacement of the domain. 
Let $\Vel:\Fluid_t \times [0,T] \to \Real^d$ and $\Pres:\Fluid_t\times [0,T] \to \Real$ 
the fluid velocity and the hydrostatic pressure. Let $\Density_f$ and $\Viscosity$ 
be the fluid density and dynamic viscosity, constant for the considered fluid model.
Finally, let $\sigma \in \mathcal{M}_d(\Real)$ be the fluid stress tensor and $g \in \Real^d$ the gravity acceleration.
The \ALE\ formulation of the Navier-Stokes equations partially decouples the geometric 
evolution of the fluid domain from that of the fluid continuum. Due to the change 
of frame, the \ALE\ time derivative substitutes the Eulerian one $\partial_t \Vel$ as 
$\partial_t \Vel = \partial_t \Vel |_\Alemap - (\partial_t x \cdot \nabla )\Vel = \partial_t \Vel |_x - ( \Vel_\Alemap  \cdot \nabla )\Vel$ 
in the momentum equation. The first term of the \ALE\ derivative corresponds 
to the time variation of the fluid velocity as seen in the arbitrary frame, 
while the second contains the relative velocity between the fluid continuum and the new reference frame.
The Navier-Stokes equations in the \ALE\ frame are:
\begin{equation}
	\begin{aligned}
		\Density_f\partial_t \Vel |_\Alemap+ \Density_f \Big((\Vel - \Vel_\Alemap)\cdot \nabla\Big) \Vel &= -\nabla \cdot \sigma + \Density_f g, \qquad &\text{in $\Fluid_t$,}\\
		\nabla \cdot \Vel &=0, \qquad &\text{in $\Fluid_t$}.
	\end{aligned}
	\label{eq:fluid}
\end{equation}

\paragraph*{The rigid body equations} 
A rigid body moving in a fluid is described by the motion of its center of mass,
given by the Newton equation, and the rotation matrix between its local frame and the laboratory 
frame, defined by the Euler equation. In this paper, the motion of a body is caused by 
fluid stresses and gravity, contact and lubrication forces that act as external 
forces $F_e$ and torques $T_e$. Let $\Solid \subset \Real^d$ be the domain occupied 
by the body, $\Density_{\Solid} \in \Real_{>0}$ 
its density and $\mass = \int_\Solid \Density_{\Solid}$ its mass. Let $\tvel: [0,T] \to \Real^d$ 
and $\angvel:[0,T]\to \Real^{d^*}$, where $d^*=1$ if $d=2$ or $d^*=3$ if $d=3$, be the 
linear and angular velocity of the rigid body as seen from the laboratory frame. 
Let $\CenterMass = \mass^{-1} \int_\Solid \Density_{\Solid} x$ be its center of mass, 
$\CenterMass\in\Real^d$, and $\Inertia= \int_\Solid \Density_{\Solid} (x-\CenterMass) \otimes (x-\CenterMass)$ 
its inertia tensor, $\Inertia\in S^{d^*}_{++}$ positive definite and symmetric.
We chose to describe the rotation matrix $\Rmat(\Angle):\Theta\to SO(d)$ using 
Euler angles $\Angle \in \Theta$, where $\Theta = [-\pi,\pi]$ if $d=2$, 
or $\Theta=[-\pi,\pi]\times[0,\pi]\times[0,\pi/2]$ if $d=3$. The 
Newton and Euler equations, describing the dynamics of a three-dimensional rigid body, are:
\begin{equation}
	\begin{aligned}
		\mass \frac{d}{dt}\tvel &= F_e-\int_{\partial \Solid} -\Pres \normal + \Viscosity(\nabla \Vel + \nabla \Vel^T)\normal,\\
		\frac{d}{dt}(\Rmat \Inertia \Rmat^T \angvel) &= T_e-\int_{\partial \Solid} [-\Pres \normal + \Viscosity(\nabla \Vel + \nabla \Vel^T) \normal ]\times (x-\CenterMass), 
	\end{aligned}
	\label{Eq:RB}
\end{equation}
where $\normal$ is the unit outward normal to $\partial \Solid$ and one has:
\begin{equation*}
	\begin{aligned}
		\frac{d}{dt} \Angle_i &= \angvel_i, \quad \text{for $i \in \{x,y,z\}$},\\
		\Rmat &=\Rmat_z(\Angle_z)\Rmat_y(\Angle_y)\Rmat_x(\Angle_x),
	\end{aligned}
\end{equation*}
where $R(\Angle_i)$ denotes the rotation matrix around axis $i \in \{x,y,z\}$ of 
angle $\Angle_i$. If $d=2$, $R(\theta)$ has the form: 
\begin{equation*}
	\begin{aligned}
		\Rmat(\Angle) &= \begin{bmatrix}
			\cos(\Angle) & \sin(\Angle)\\
			-\sin(\Angle) & \cos(\Angle)
		\end{bmatrix},
	\end{aligned}
\end{equation*}
while in three dimensions:
$$	
	\Rmat_z(\Angle_z) = \begin{bmatrix}
		\cos(\Angle_z) & -\sin(\Angle_z) & 0\\
		\sin(\Angle_z) & \cos(\Angle_z)& 0 \\
		0& 0 & 1
	\end{bmatrix},
	\Rmat_y(\Angle_y) = \begin{bmatrix}
		\cos(\Angle_y) &0 & \sin(\Angle_y) \\
		0& 1& 0\\
		-\sin(\Angle_y) & 0& \cos(\Angle_y)
	\end{bmatrix},\\
$$
$$
	\Rmat_x(\Angle_x) = \begin{bmatrix}
		1 & 0 & 0\\
		0 &\cos(\Angle_x) & -\sin(\Angle_x) \\
		0 &\sin(\Angle_x) & \cos(\Angle_x) 
	\end{bmatrix}.
$$

\paragraph*{Fluid-solid interaction} 
The interaction between the rigid body and the fluid is ensured by 
the balance of stresses by the Newton and Euler equations \eqref{Eq:RB} and a 
coupling condition at the interface, imposing the continuity of the velocities:
\begin{equation}
	\Vel = \tvel + \angvel \times (x-\CenterMass) \quad \text{on $\partial \Fluid_t \cap \partial\Solid$}.
\end{equation}

Finally, the definition of the $\Alemap_t(X)$ requires the computation of the 
displacement $x(t,X) $ in 
the reference domain. In order to find $x(t,X)$, we solve:
\begin{equation}
	\begin{aligned}
		\nabla \cdot([1+\tau(X)]\nabla_X x(t,X)) &= 0 \quad &\text{in $\Fluid_0$}, \\
		x(t,X) &= g(t,X) \quad &\text{in $\partial \Fluid_0$}, 
	\end{aligned}
	\label{Eq:ALE}
\end{equation}
where $g(t,X)  = \int_0^t \tvel + \angvel \times(X-X^{CM})$ is the rigid displacement 
of $\Solid$ and $\tau(X)$ is a space-dependent coefficient, related to the 
volume of the simplexes in the triangulation. 

\paragraph*{Numerical solution}
The problem is discretized over a triangulation of the fluid domain. The discrete 
fluid velocity $\Vel_h$ and pressure $\Pres_h$ belong to the inf-sup stable Taylor-Hood 
space, where $P_2$ continuous finite elements are chosen for the velocity and $P_1$ 
continuous finite elements for the pressure. The discrete \ALE\ map $\Alemap_h^t$ 
is also discretized using $P_1$ continuous finite elements.
The coupling condition of the fluid-solid interaction is encoded in the finite 
element test spaces, as proposed in \cite{maury_direct_1999}, which 
leads to the construction of an operator $P$ that couples the velocities at the 
interface, at the discrete level.
In our computations, the interface is conforming and the solution of equation 
\eqref{Eq:ALE} handles small mesh deformations, with $\tau(X)$ a piecewise 
constant coefficient, defined on each element $e$ of the mesh 
as $ \tau\big |_e=(1-V_{min}/V_{max})/(V_e/V_{max})$, where $V_{max}$, $V_{min}$  
and $V_e$ are the volumes of the largest, smallest and current element of the domain 
discretization \cite{kanchi_3d_2007}. 
In case of larger domain deformations, we perform remeshing and preserve the 
discretization of the interface at the same time.\\

\section{Collision model} \label{collisionmodel}

In this section, we will detail our collision model characterized by two phases: a 
collision detection algorithm and a lubrication model. 
We do not add a scheme defining direct contact between body surfaces, since our lubrication 
model avoids this kind of interaction. We will first describe the three types of 
rigid bodies for which we have 
developed our collision model. Then, we will detail the collision detection algorithm 
before defining the lubrication force. 

\subsection{Rigid body types}

We adjust the collision model for three different types of rigid bodies. 
First, we consider spherical bodies in two and three dimensions. As collision models 
for spherical bodies are common in literature, we use this case for validation. 
Then, we adapt this first model for bodies of complex shape. In this paper, 
ellipses are used for the two-dimensional simulations and ellipsoids for simulations 
in three dimensions. At last, we consider articulated bodies, and in particular the 
three-sphere swimmer \cite{najafi_simplest_2004}. This type of swimmer is composed 
of three spheres of the same size that are connected by rods. The swimmer extends and 
retracts these rods to move. First the left sphere is retracted, then the right sphere. 
Finally, the left sphere is extended before extending the right sphere. This sequence 
of four movements result in a straight motion.  
	
\subsection{Collision detection}

The first part of our model is a collision detection algorithm. Collision detection 
is an important task since it allows identifying the pairs of bodies that are actually 
interacting. Two bodies are interacting when the distance between their surfaces is 
smaller than the width of the collision zone, denoted $\rho$. In absence of detection phase, 
collision forces are computed for each pair 
of bodies present in the domain, independently of whether there will be an interaction. 
Most collision detection algorithms include two phases. First, the broad-phase identifies 
pairs of bodies in close neighborhood. Then, the distance between the identified 
pairs is computed explicitly during the narrow-phase. According to this distance, 
one concludes if collision forces are applied. Our collision detection algorithm 
consists only of the narrow-phase. Thus, the distance 
between each pair of bodies is computed. Depending on the body shape, this is 
done in two different ways. We will first describe it for spherical and then for 
complex shaped or articulated bodies. Both methods require some parameters 
defined in a pre-processing step, which will also be described. 

\begin{algorithm}
	\footnotesize
    \caption{Pre-process phase}
    \begin{algorithmic}
	\footnotesize

	\State \textbf{Input:} Computational domain: mesh $\mbox{M}$ 
	\State \textbf{Output:} Set of bodies identifiers: $\mbox{bodyIds}$, set of centers of mass: $\mbox{massCenters}$, 
	set of bodies markers: $\mbox{bodyMarkers}$, set of radii: $\mbox{radii}$, set of centers of mass of imaginary bodies: $\mbox{ImagmassCenters}$, and the fluid marker: $\mbox{fluidMarker}$
	
	\smallskip

	\For{$\Solid$ in M.bodies}
		\State $\mbox{bodyIds}$.append($\Solid$.id())  
		\State $\mbox{massCenters}$.append($\Solid$.massCenter())
		\State $\mbox{bodyMarkers}$.append($\Solid$.boundary().name())
		
		\medskip
		
		\If{type = spherical} 
        	\State $\mbox{radii}$.append($\Solid$.radius()) 
			\State $\mbox{ImagmassCenters}$.append($\Solid$.ImagmassCenters()) 
		\EndIf

		\smallskip

    \EndFor

	\smallskip

	\State $\mbox{fluidMarker} \leftarrow \Fluid_t$.boundary().name()  
    \end{algorithmic}
    \label{algo_preprocess}
    
\end{algorithm}
\normalsize

\paragraph*{Spherical bodies}
The distance $d_{ij}$ between the surfaces of two spherical bodies $\Solid_i$ and $\Solid_j$, $i \ne j$, 
is determined by first computing the distance between their mass centers $\CenterMass_i$ and $\CenterMass_j$ 
and then subtracting their radii $r_i$ and $r_j$:
$$
d_{ij} = || \CenterMass_{j} - \CenterMass_i ||_2 - r_i - r_j.
$$
When this distance is smaller than a given parameter $\rho$, representing the width 
of the collision zone (detailed in subsection \ref{lubmodel}), then the identifiers $i$, $j$ 
and the distance $d_{ij}$ are stored in a collision map. 
Stored data is used during the phase where lubrication forces are computed.\\
The same formula can be used to get the distance $d_i$ between the surface of one 
body $\Solid_i$ and the boundary of the fluid domain $\partial \Fluid_t$:
$$
d_{i} = || \CenterMass_{i' } - \CenterMass_{i}  ||_2 - 2r_i.
$$
where $\CenterMass_{i'}$ is the center of mass of the nearest imaginary body placed on the outside 
of the fluid boundary \cite{wan_direct_2006}. This imaginary body has same radius as the real body. Again, if $d_i$ 
is smaller than the width $\rho$, then the data $i$ and $d_i$ are stored. 
The parameters used in these two equations, i.e. body identifiers, mass centers and radii, 
are defined in a pre-processing step given by algorithm \ref{algo_preprocess}. 

\begin{algorithm}
	\footnotesize
    \caption{Collision detection algorithm for spherical bodies}
    \begin{algorithmic}
	\footnotesize
	\State \textbf{Input:} Output of pre-process phase
	\State \textbf{Output:} Collision map
	
	\smallskip 
	
	\For{$i$ in  bodyIds}

		\smallskip 

		\For{$j$ in bodyIds, $j > i$}
			\State $d_{ij} \leftarrow ||\mbox{massCenters}_{i} - \mbox{massCenters}_{j}||_2 - \mbox{radii}_i - \mbox{radii}_j$
			\If{$d_{ij} \leq \rho$}
				\State Store $i$, $j$ and $d_{ij}$ in collision map 
			\EndIf
		\EndFor	

		\smallskip 

		\For{ImagCoord in $\mbox{ImagmassCenters}_{i}$}
			\State $d_{i} \leftarrow ||\mbox{massCenters}_{i} - \mbox{ImagCoord}||_2 - 2\mbox{radii}_i$
			\If{$d_{i} \leq \rho$}
				\State $j \leftarrow -1$
				\State Store $i$, $j$, ImagCoord and $d_{i}$ in collision map 
			\EndIf
		\EndFor	

		\smallskip 

    \EndFor	
    \end{algorithmic}
    \label{algo_detection_1}
\end{algorithm}
\normalsize

\paragraph*{Complex shaped or articulated bodies}
When considering complex shaped or articulated bodies, then no explicit expression 
is provided to compute the distance between the surfaces. 
To deal with this issue, we use the fast marching method, based on a level 
set algorithm introduced by J.A Sethian \cite{sethian_fast_1996}. By applying it to a body $\Solid_i$, it returns 
the distance field $D_i$ from that body to the rest of the domain. Thus, the distance 
field $D_i$ is equal to zero at the body's boundary and has a positive value elsewhere. 
Given the distance fields $D_i$ and $D_j$, then the minimum distance between the 
surfaces $\partial \Solid_i$ and $\partial \Solid_j$ of the corresponding bodies is found by:
$$
d_{ij} = ||\arg\min_{x \in \partial \Solid_j}  D_i(x) - \arg\min_{x \in \partial \Solid_i} D_j(x)||_2. 
$$
The minima $X_i = \arg\min_{x \in \partial \Solid_i} D_j(x)$ and $X_j = \arg\min_{x \in \partial \Solid_j}  D_i(x)$ 
represent the contact points of the bodies $\Solid_i$ and $\Solid_j$, i.e. the coordinates 
of the surface points where the bodies will interact. To get the minimal distance between a 
body surface and the domain boundary, we use: 
$$
d_{i} = ||\arg\min_{x \in \partial \Fluid_t} D_i(x) - \arg\min_{x \in \partial \Solid_i} D_{\Fluid_t}(x)||_2 \mbox{,}
$$
where $D_{\Fluid_t}$ the distance field obtained by applying the fast marching method 
to the domain boundary $\partial \Fluid_t \backslash ( \cup_i \partial \Solid_i ) $. 
The parameters that will be stored if $d_{ij} \leq \rho$ or $d_i \leq \rho$, 
are the identifiers, the minimal distance and the contact points. \newline 

\begin{algorithm}
	\footnotesize
    \caption{Collision detection algorithm for complex shaped bodies}
    \begin{algorithmic}
	\footnotesize

	\State \textbf{Input:} Output of pre-process phase
	\State \textbf{Output:} Collision map
	
	\smallskip 

	\For{$i$ in  bodyIds}
		\State $D_i \leftarrow \mbox{narrow-band}(\mbox{bodyMarkers}_{i}, d_{max} = 1.5 \rho)$ 
    \EndFor

	\State $D_{\Fluid_t} \leftarrow \mbox{narrow-band}(\mbox{fluidMarker}, d_{max} = 1.5 \rho)$ 

	\smallskip

	\For{$i$ in  bodyIds}
		
		\smallskip

		\For{$j$ in bodyIds, $j > i$}
			\State $d_{ij} \leftarrow ||X_j - X_i||_2$
			\If{$d_{ij} \leq \rho$}
				\State Store $i$, $j$, $d_{ij}$, $X_i$ and $X_j$ in collision map 
			\EndIf
		\EndFor	

		\smallskip 

		\State $d_{i} \leftarrow  ||X_{\Fluid_t} - X_i||_2$
		\If{$d_{i} \leq \rho$}
			\State $j \leftarrow -1$ 
			\State Store $i$, $j$, $d_{i}$, $X_i$ and $X_{\Fluid_t}$ in collision map 
		\EndIf
    \EndFor
	
    \end{algorithmic}
    \label{algo_detection_2}
\end{algorithm}
\normalsize

When the fluid domain is large, the fast marching algorithm gets computationally expensive, 
especially in three dimensions. To accelerate the method, we develop the narrow-band approach, 
which computes the distance field only in a near neighborhood of the body boundary. 
This neighborhood is set to a predefined threshold $d_{max}$. As collision forces 
are only applied on a collision zone of width $\rho$, the threshold $d_{max}$ can 
be defined as function of $\rho$. Once the threshold distance $d_{max}$ is reached, the narrow-band 
approach assigns a default value $\delta$ to the distance field, 
corresponding to the maximum value reached:
\begin{center}
	\begin{equation*}
		D_i^{NB} = \left\{\begin{array}{rcl}
		0 \;, \quad &\mbox{on }& \partial \Solid_i , \\
		D_i \;, \quad &\mbox{for }& D_i \leq d_{max} ,\\
		\delta\;, \quad &\mbox{elsewhere}&   .
		\end{array}\right.\;
	\end{equation*}
\end{center}

The narrow-band algorithm can be executed in parallel. A performance study of this 
approach is performed in section \ref{experiments}.
In our implementation, the narrow-band fast marching function takes as parameter 
the boundary marker of the concerned body. Therefore, all boundary markers are 
stored during the pre-processing phase, presented in Algorithm \ref{algo_preprocess}. \\ 

When considering articulated bodies, the same formulas as for the complex shaped 
case are applied to each component of the body. Thus, when considering a pair of 
three-sphere swimmers, the distances between each sphere of one swimmer and all 
those of the second are computed. All steps of the collision detection algorithm 
for both methods are given by Algorithm \ref{algo_detection_1} for spherical and 
by Algorithm \ref{algo_detection_2} for complex shaped bodies. The computation of the  imaginary 
bodies centers is only possible when the dimensions of the fluid domain, which 
must be rectangular, are known. For this reason, we have developed a more generic 
algorithm for spherical bodies, where the distance between a body and the fluid boundary 
is determined using the narrow-band approach. The initial algorithm can still be used 
for simulations of spherical bodies in rectangular domains, because its execution costs are lower.  
	
\subsection{Collision force} \label{lubmodel}

The second part of our collision model consists in the computation of collision forces. The lubrication model is based 
on a short-range repulsive force, a scheme introduced by R. Glowinski \cite{glowinski_linite_nodate}. 
The force is activated when the distance between two body surfaces is less than 
the parameter $\rho$. This parameter represents the width of the collision zone, i.e. 
the zone where collision forces must be applied to prevent the bodies from overlapping and 
to avoid direct contact. 
The definition of the force is quite simple: it is parallel to the vector that connects 
the contact points of the bodies and its intensity increases as the distance decreases. 
We will first give the definition of this force for the simple case of spherical bodies. 
Then, we detail the required modifications, in particular the addition of its 
torque, for complex shaped and articulated bodies. \newline

\paragraph*{Spherical bodies}
For the repulsive force definition in the case of spherical bodies, we rely on articles \cite{wang_drafting_2014, usman_analysis_2018, wan_direct_2006}. 
The force to model the interaction between two bodies $\Solid_i$ and $\Solid_j$ is given by:
$$
\overrightarrow{F_{ij}} = \frac{1}{\epsilon} \overrightarrow{\CenterMass_j \CenterMass_i} (\rho - d_{ij})^2 \mathds{1}_{d_{ij} \leq \rho}.
$$
In the same way, the interaction between a body $\Solid_i$ and the domain boundary $ \partial \Fluid_t$ is defined by:
$$
\overrightarrow{F_{i \Fluid_t}} = \frac{1}{\epsilon_{\Fluid}} \overrightarrow{\CenterMass_{i'} \CenterMass_i} (\rho - d_{i})^2 \mathds{1}_{d_{i} \leq \rho}.
$$
Both equations contain a quadratic activation term. The vector connecting the mass 
centers gives the direction of the force, and the stiffness parameters $\epsilon$ and $\epsilon_{\Fluid}$ 
describe its intensity for body-body and body-domain interactions.
Finding the optimal values for these two parameters is not trivial, 
since their values depend on fluid and body properties. We refer to the indications 
given by article \cite{wan_direct_2006}: when the width of the collision zone is fixed to $\rho = 0.5h \sim 2.5h$, 
where $h$ is the mesh size, and the ratio between the body and fluid density is 
equal to $1$, then one can suppose $\epsilon \approx h^2$ and 
$\epsilon_{\Fluid} \approx \frac{\epsilon}{2}$.
As already mentioned, these forces are only computed for pairs stored during 
collision detection phase, i.e. for pairs whose distance is less than $\rho$, 
in which cases the indicator functions $\mathds{1}_{d_{ij} \leq \rho} \equiv 1$ and $\mathds{1}_{d_{i} \leq \rho} \equiv 1$. 
The total repulsion force applied on one 
body $\Solid_i$ is defined by:
$$
\overrightarrow{F_i} = \sum_{(i,j)|d_{ij} \leq \rho} \overrightarrow{F_{ij}} +  \sum_{i | d_i \leq \rho} \overrightarrow{F_{i\Fluid_t}} ,
$$
for all pairs $\Solid_i$ - $\Solid_j$
and $\Solid_i$ - $\partial \Fluid_t$ present in the collision map after detection algorithm.
The total repulsion force $\overrightarrow{F_i}$ is added to the external forces 
$\overrightarrow{F_e}$ of the Newton equation in \eqref{Eq:RB}, describing the linear velocity of the body $\Solid_i$.  

\begin{algorithm}
	\footnotesize
    \caption{Collision force algorithm for spherical bodies}
    \begin{algorithmic}
		
		\State \textbf{Input:} Output of pre-process phase and collision map
		\State \textbf{Output:} Set of collision forces: $\mbox{F}$, and set of collision torques: $\mbox{T}$

		\smallskip
	
		\For{pair in  collision map}

			\If{pair.$j = -1 $ } 
				\State $F_{i \Fluid_t} \leftarrow \frac{1}{\epsilon_{\Fluid}} (\mbox{massCenters}_{\mbox{pair}.i} - \mbox{pair.ImagCoord}) (\rho - \mbox{pair.}d_{i})^2 $
				\State $F_{\mbox{pair}.i}$.add($F_{i \Fluid_t}$) 
			\Else 
				\State $F_{ij} \leftarrow \frac{1}{\epsilon} (\mbox{massCenters}_{\mbox{pair}.i} - \mbox{massCenters}_{\mbox{pair.}j}) (\rho - \mbox{pair.}d_{ij})^2 $
				\State $F_{\mbox{pair}.i}$.add($F_{ij}$),  $F_{\mbox{pair}.j}$.add($- F_{ij}$) 
			\EndIf
		
    	\EndFor
    \end{algorithmic}
    \label{algo_force_1}
\end{algorithm}

\paragraph*{Complex shaped or articulated bodies}

For the case of non-spherical bodies, the direction of the vector connecting the mass 
centers $\CenterMass_i$ and $\CenterMass_j$ no longer corresponds to that of the vector connecting the 
contact points $X_i$ and $X_j$. Thus, the equations are 
slightly modified: 

$$
\vec{F_{ij}} = \frac{1}{\epsilon} \overrightarrow{X_j X_i } (\rho - d_{ij})^2 \mathds{1}_{d_{ij} \leq \rho}, \qquad
\vec{F_{i \Fluid_t}} = \frac{1}{\epsilon_{\Fluid}} \overrightarrow{X_{\Fluid_t} X_i} (\rho - d_{i})^2 \mathds{1}_{d_{i} \leq \rho}.
$$
The values of the parameters $\epsilon$ and $\epsilon_{\Fluid}$ have to be chosen smaller than 
those of the previous system. These modified expressions imply that the collision forces will lead to body rotations. 
The torque to model these rotations is added to the external torques $T_e$ of the 
Euler equation in \eqref{Eq:RB}, and is defined by: 
$$
\overrightarrow{T_i} = - \overrightarrow{\CenterMass_i X_i} \times \overrightarrow{F_i}.
$$
This torque is computed for each body present in the collision map. 
The same approach is used for articulated bodies. The forces and torques are 
added for each element of the collision map.
The algorithms given by Algorithm \ref{algo_force_1} and \ref{algo_force_2} illustrate the implementation of the 
lubrication model for spherical and complex shaped bodies.

\begin{algorithm}
	\footnotesize
    \caption{Collision force algorithm for complex shaped bodies}
    \begin{algorithmic}

		\State \textbf{Input:} Output of pre-process phase and collision map
		\State \textbf{Output:} Set of collision forces $\mbox{F}$, and set of collision torques $\mbox{T}$

		\smallskip
	
		\For{pair in  collision map}
			\If{pair.$j = -1 $ } 
				\State $F_{i \Fluid_t} \leftarrow \frac{1}{\epsilon_{\Fluid}} (\mbox{pair}.X_{i} - \mbox{pair.}X_{\Fluid_t}) (\rho - \mbox{pair.}d_{i})^2 $
				\State $F_{\mbox{pair}.i}$.add($F_{i \Fluid_t}$) 
				\State $T_{\mbox{pair}.i}$.add($- (\mbox{pair}.X_{i} - \mbox{massCenters}_{\mbox{pair}.i}).cross(F_{i \Fluid_t})$)  
			\Else 
				\State $F_{ij} \leftarrow \frac{1}{\epsilon} (\mbox{pair}.X_{i} - \mbox{pair.}X_{j}) (\rho - \mbox{pair.}d_{ij})^2 $
				\State $F_{\mbox{pair}.i}$.add($F_{ij}$),  $F_{\mbox{pair}.j}$.add($- F_{ij}$) 
				\State $T_{\mbox{pair}.i}$.add($- (\mbox{pair}.X_{i} - \mbox{massCenters}_{\mbox{pair}.i}).cross(F_{ij})$) 
				\State $T_{\mbox{pair}.j}$.add($- (\mbox{pair}.X_{j} - \mbox{massCenters}_{\mbox{pair}.j}).cross(-F_{ij})$) 
			\EndIf
	
		\EndFor	
    \end{algorithmic}
    \label{algo_force_2}
\end{algorithm}

\section{Numerical experiments} \label{experiments}

\subsection{Performance of narrow-band approach}

The following test illustrates the performance of the narrow-band approach of the 
fast marching method, in two and three dimensions. The geometry consists in two 
spheres of same center but different radii. 
The radius of the inner and outer sphere is respectively set to $r_i=0.1$ and $r_o=2.0$. 
First, we use the fast marching method to determine the distance field from the inner 
sphere boundary to the outer sphere boundary. This is equivalent to perform the narrow-band 
approach on the inner sphere by setting the threshold distance $d_{max}$ to $d_{max} = r_o - r_i = 1.9$. 
Then, we continue to apply the narrow-band approach on the inner sphere but considering 
smaller thresholds. For the two-dimensional test the thresholds are set to 
$d_{max} = \{ 0.0625, 0.125, 0.25, 0.5, 1.0\}$, and we fix the mesh size to   
$h = 0.01$. The mesh contains $E_{2D} = 294918$ elements. In three dimensions, the mesh 
has $E_{3D} = 1385385$ elements, using a size equal to $h = 0.05$. The set of thresholds 
is  $d_{max} = \{0.125, 0.25, 0.5, 1.0\}$.\newline

The table \ref{Fig:Performance} shows the execution time of the narrow-band 
approach for the different thresholds as well as the speed-up in two and three dimensions. 
We add the number of band elements $E_b$, i.e. the number of mesh elements in the zone 
where the distance field is actually computed, and the element ratio 
$\frac{E_{nD}}{E_b}$, where $n$ the dimension.
It can be observed that the execution time decreases significantly when considering 
smaller thresholds. This test illustrates the importance of the narrow-band approach for 
our collision model. 

\begin{table}[h!]
	
	\begin{center}
		
		\begin{tabular}{|r|r|r|r|r|} 
		    \hline  
		    Threshold $d_{max}$ & Band elements $E_b$& Element ratio $\frac{E_{2D}}{E_b}$ & Execution time in $2$D & Speedup \\  
		    \hline 
		    $1.9$ & $294918$ &  $1.00$ & $1.959$ s &   /  \\  
		    \hline
		    $1.0$ & $89898$ & $3.28$ & $0.722$ s &  $2.71$  \\ 
		    \hline 
		    $0.5$ & $27082$ & $10.89$ & $0.411$ s &  $4.76$  \\ 
		    \hline 
		    $0.25$ & $9430$ & $31.27$ & $0.316$ s &  $6.20$  \\ 
		    \hline 
		    $0.125$ & $4104$ & $71.86$ & $0.276$ s & $7.09$  \\ 
			\hline 
			$0.0625$ & $2262$ & $130.38$ & $0.261$ s & $7.50$  \\ 
		    \hline 
		    Threshold $d_{max}$ & Band elements $E_b$ & Element ratio $\frac{E_{3D}}{E_b}$& Execution time in $3$D & Speedup \\   
		    \hline 
		    $1.9$ & $1385385$ &  $1.00$ & $22.851$ s   & /  \\  
		    \hline
		    $1.0$ & $254365$ & $5.45$ & $5.869$ s & $3.89$  \\ 
		    \hline 
		    $0.5$ & $43784$ &  $31.64$ & $2.368$ s & $9.65$  \\ 
		    \hline 
		    $0.25$ & $10238$ & $135.32$ & $1.815$ s & $12.59$  \\ 
		    \hline 
		    $0.125$ & $4537$ &  $305.35$ & $1.759$ s & $12.99$  \\ 
		    \hline 
	        \end{tabular}
	
    \end{center}

	\caption{Results of performance test for the narrow-band 
	approach in two and three dimensions. The table shows the execution time and 
	speedup for different threshold values $d_{max}$.}

	\label{Fig:Performance}
\end{table}

\subsection{Performance of the parallel collision algorithm}

The aim of this second numerical experiment is to observe the performance of 
the parallel implementation of the collision algorithm using the narrow-band fast 
marching method as well as to analyze the execution time proportion attributed to 
the collision model during the fluid-solid interaction resolution. 
The choice of the preconditioners, a direct LU solver in two dimensions and a 
block preconditioner in $3$D, are detailed in \cite{berti_fluid-rigid_nodate}.
We consider different numbers 
of spherical bodies distributed in a rectangular domain filled with a steady fluid. 
No gravity is applied and the 
magnitude of the collision forces is chosen close to zero, such that the bodies do 
not move during the simulation. 
Each domain configuration is simulated for ten time
iterations in sequential and in parallel, using different numbers of processors. 
For each case various information 
are displayed: the number of bodies, the number of mesh nodes, the number of body-body 
and body-wall interactions during each iteration as well as the speed-up 
of the collision algorithm (pre-process phase, collision detection algorithm, and 
computation of collision force) 
after the ten iterations. This speed-up $\frac{T_{np1}}{T_{npN}}$ is given by the 
ratio between the sequential execution time and the parallel execution time, 
when running the simulation on $N$ processors. 
The results are given in tables \ref{per2d} 
and \ref{per3d}.

\begin{table}[h!]
	
	\begin{center}
		\begin{tabular}{|r|r|r|r|r|r|r|r|} 
		    \hline  
		    Bodies & Mesh nodes & Body-body & Body-wall & $\frac{T_{np1}}{T_{np4}}$  & $\frac{T_{np1}}{T_{np16}}$ & $\frac{T_{np1}}{T_{np24}}$ \\  
		    \hline 
		    $1$ & $16462$ &  $0$ & $1$ &   $3.14$  & $8.22$ & $11.04$ \\  
		    \hline
		    $25$ & $17780$ & $0$ & $16$ &  $1.88$   & $2.06$ & $2.12$ \\ 
		    \hline 
		    $49$ & $18284$ & $28$ & $24$ &  $1.90$   & $2.08$ & $2.15$ \\ 
		    \hline 
		    $81$ & $18813$ & $72$ & $32$ & $1.92$   & $2.09$ & $2.13$ \\ 
			\hline 
			$100$ & $19139$ & $100$ & $36$  & $1.81$   & $2.01$ & $2.10$ \\ 
			\hline
	        \end{tabular}
	
    \end{center}

	\caption{Results of the parallelization performance of the collision algorithm 
	in two dimensions. The table shows the speed-up of the execution time when 
	running the simulations in parallel.}

	\label{per2d}
\end{table} 

\begin{table}[h!]
	
	\begin{center}
		\begin{tabular}{|r|r|r|r|r|r|r|r|} 
		    \hline  
		    Bodies & Mesh nodes & Body-body & Body-wall & $\frac{T_{np1}}{T_{np4}}$ & $\frac{T_{np1}}{T_{np16}}$ & $\frac{T_{np1}}{T_{np24}}$ \\  
		    \hline 
		    $1$ & $17073$ &  $0$ & $1$ &   $2.23$  & $6.49$ & $7.99$ \\  
		    \hline
		    $27$ & $18213$ & $0$ & $26$ &  $1.94$   & $3.99$ & $5.03$ \\ 
		    \hline 
		    $64$ & $18789$ & $0$ & $56$ &  $1.98$  & $3.98$ & $5.22$ \\ 
		    \hline 
		    $125$ & $20301$ & $70$ & $98$ & $1.99$   & $3.59$ & $5.13$ \\ 
			\hline
	        \end{tabular}
	
    \end{center}

	\caption{Parallelization performance in three dimensions. 
	The sequential execution time is compared to the execution time of the simulation 
	when running it in parallel using different numbers of processors.}

	\label{per3d}
\end{table} 

One can observe that the parallel execution reduces the execution times; 
the speed-up increases with an increasing number of processors. For multiple interactions, 
the speed-up remains almost constant independently of the number 
of bodies.
The execution time spent in collision detection 
represents the majority of the total execution time of the algorithm. 
For the simulations conducted in this section, the execution time of the detection 
phase represents $99\%$ of the total collision 
time. 
The difference in speed-ups between the two and three-dimensional case could be explained by the 
larger execution time taken by the distance computation, since the number of node-node 
connections in $3$D is, on average, larger than in $2$D.
Furthermore, the execution time of the two-dimensional 
collision model represents up to $20 \%$ of the total time for the fluid-solid 
interaction resolution, when considering multiple interactions and using a direct solver.
In three dimensions, using a block preconditioner, the proportion associated to the 
collision algorithm does not exceed $5\%$. 
\section{Validation} \label{validation}

To validate our model we consider the motion of a circular solid in an incompressible 
Newtonian fluid in two and three dimensions.\newline

\begin{figure}
	
	\begin{center}
	\begin{subfigure}{0.24\textwidth}
		\centering
		\includegraphics[scale=0.27]{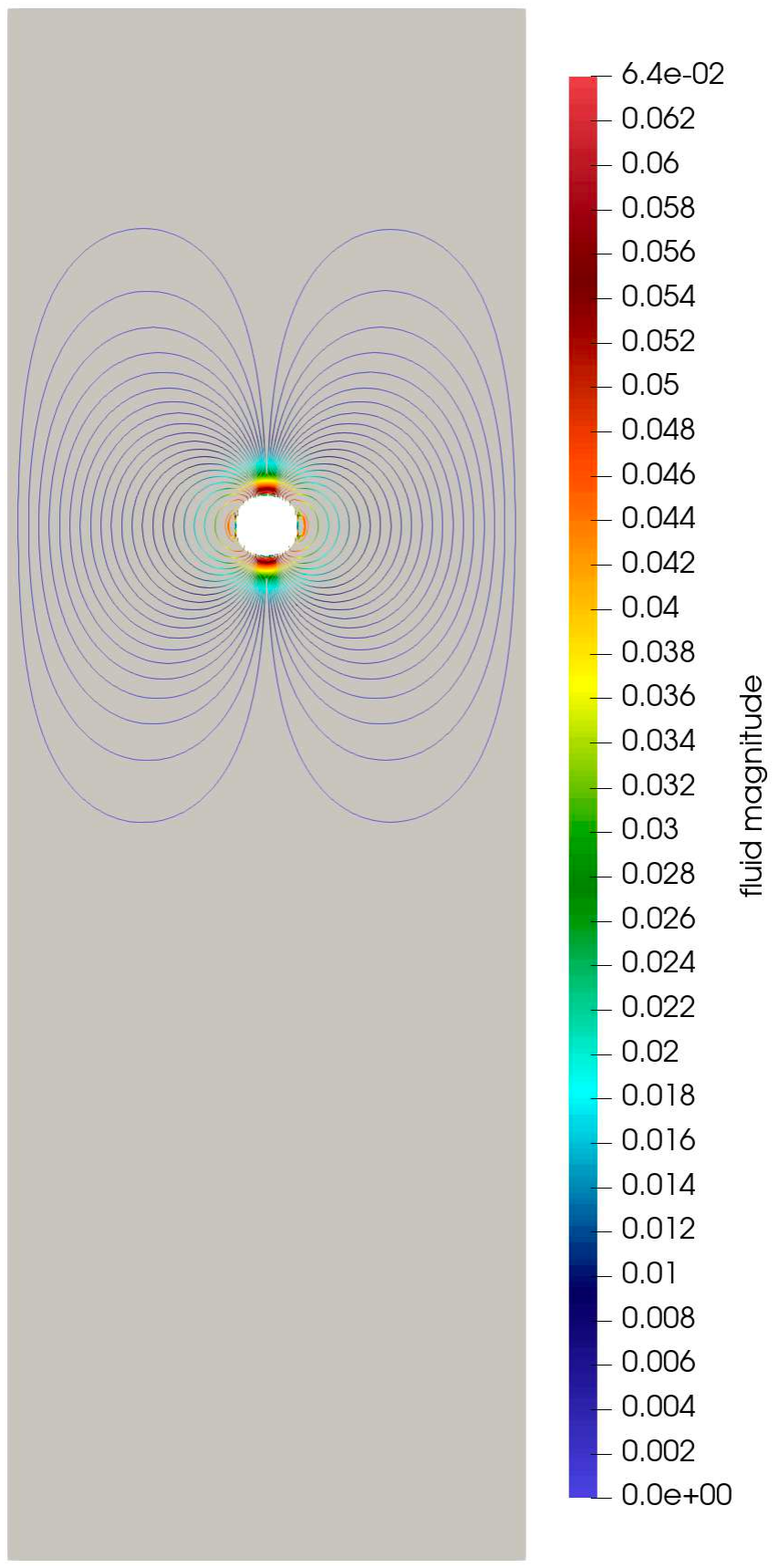}
		\caption{$t=0.0 s$}
		\label{fig:validation2d_1}
	\end{subfigure}
	\begin{subfigure}{0.24\textwidth}
		\centering
		\includegraphics[scale=0.27]{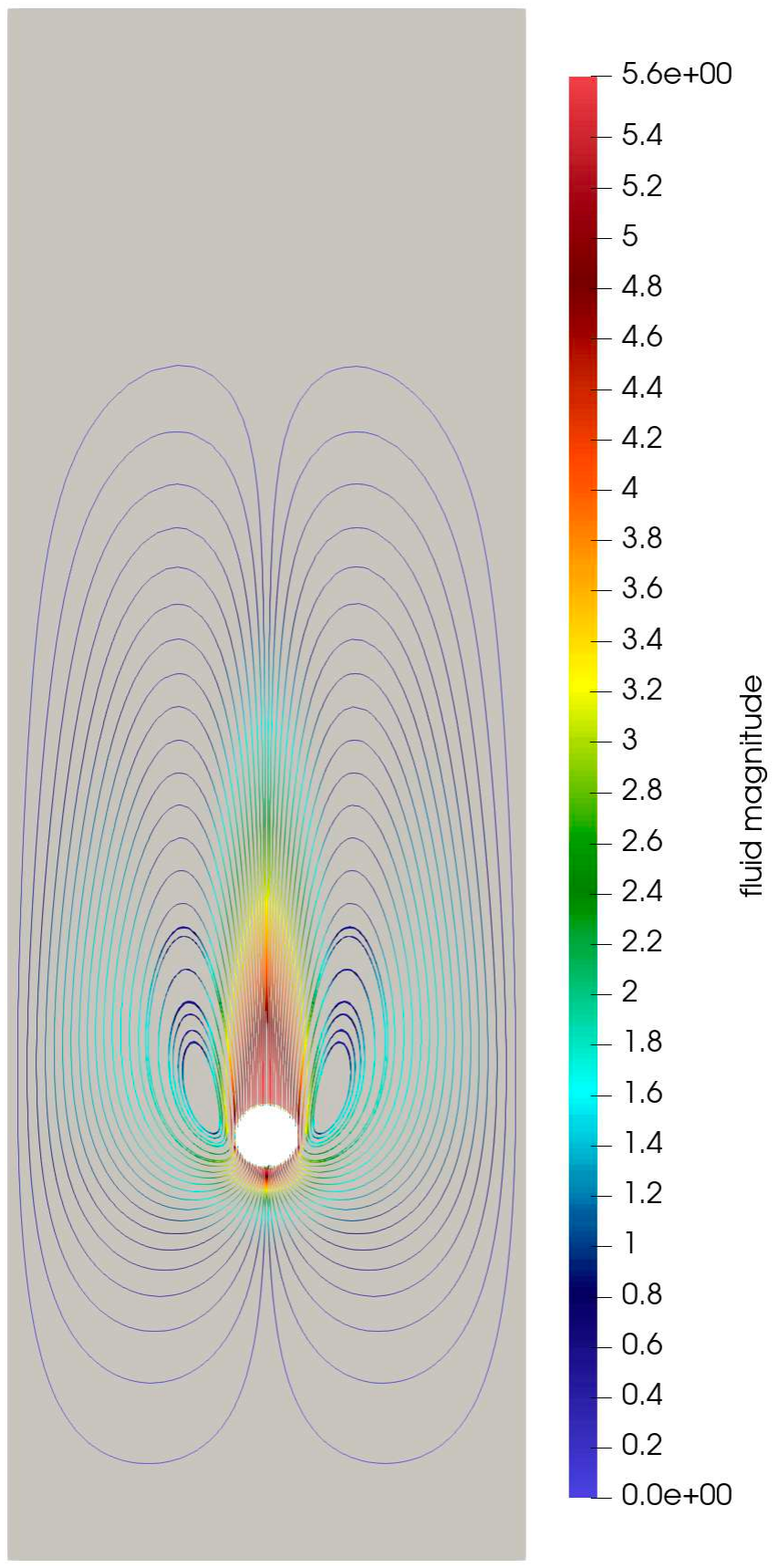}
		\caption{$t=0.5 s$}
		\label{fig:validation2d_2}
	\end{subfigure}
	\begin{subfigure}{0.24\textwidth}
		\centering
		\includegraphics[scale=0.27]{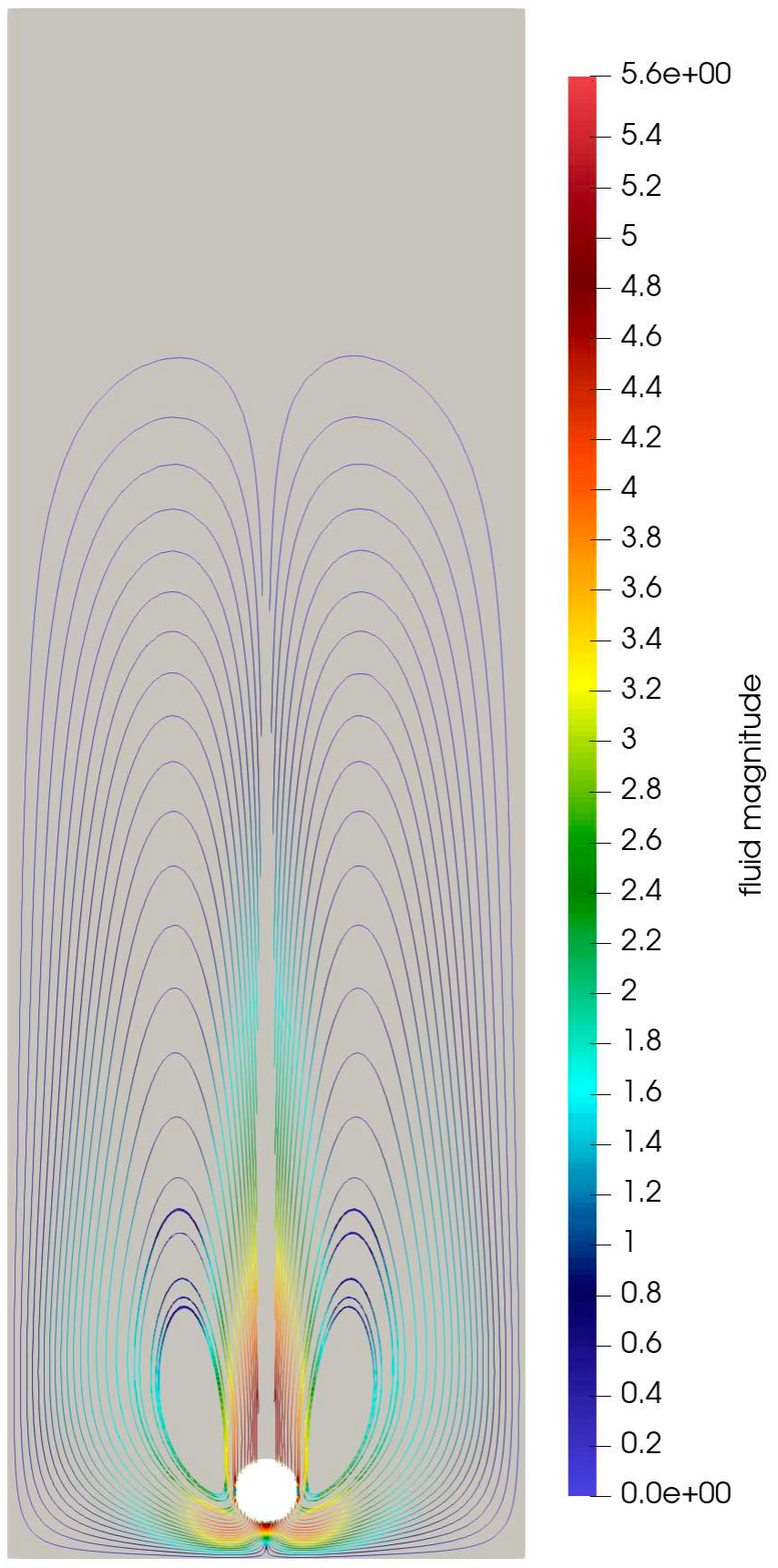}
		\caption{$t=0.75 s$}
		\label{fig:validation2d_3}
	\end{subfigure}
	\begin{subfigure}{0.24\textwidth}
		\centering
		\includegraphics[scale=0.27]{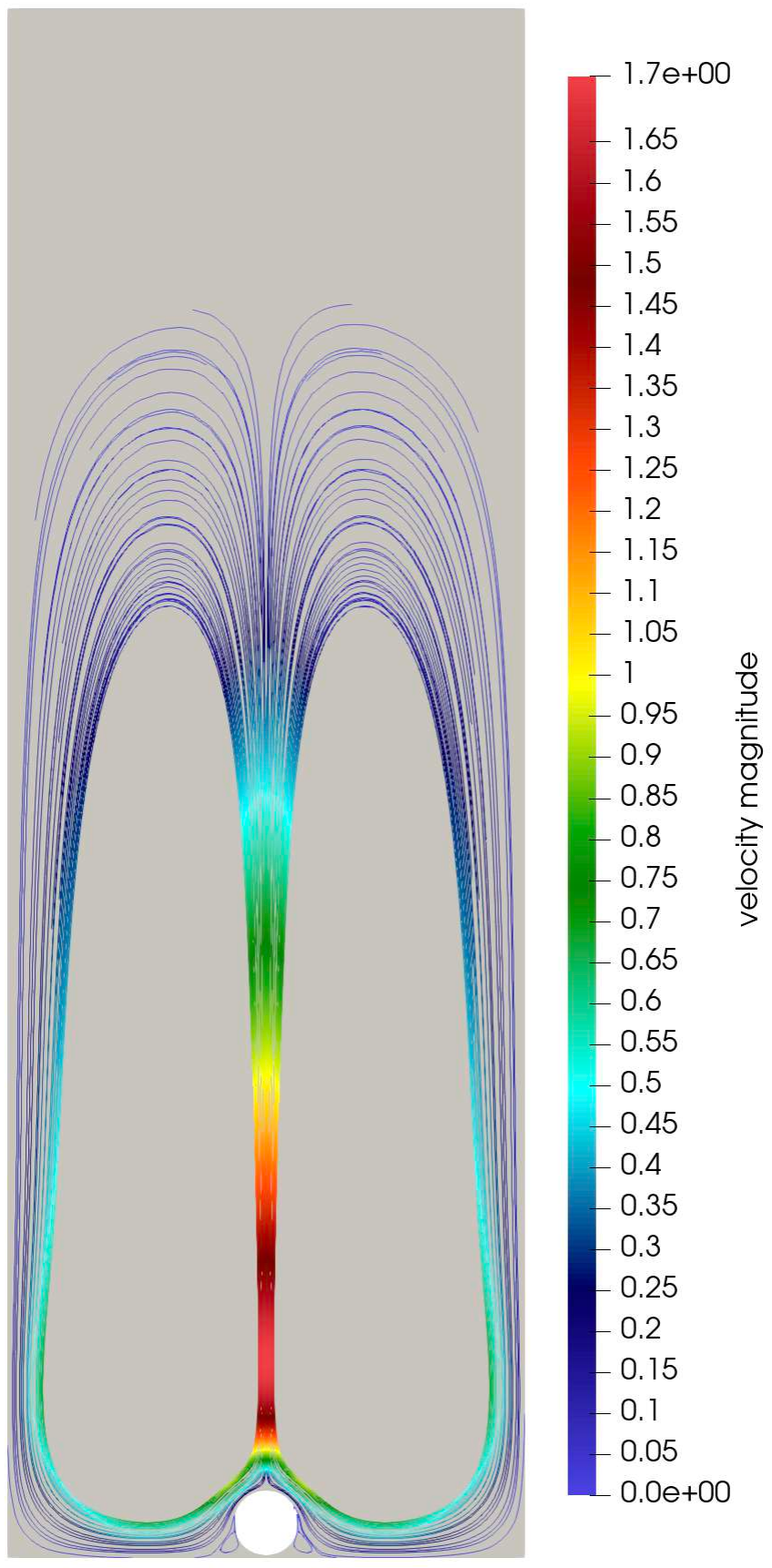}
		\caption{$t=1.0 s$}
		\label{fig:validation2d_4}
	\end{subfigure}
	\end{center}

	\centering  
	\includegraphics[scale=0.25]{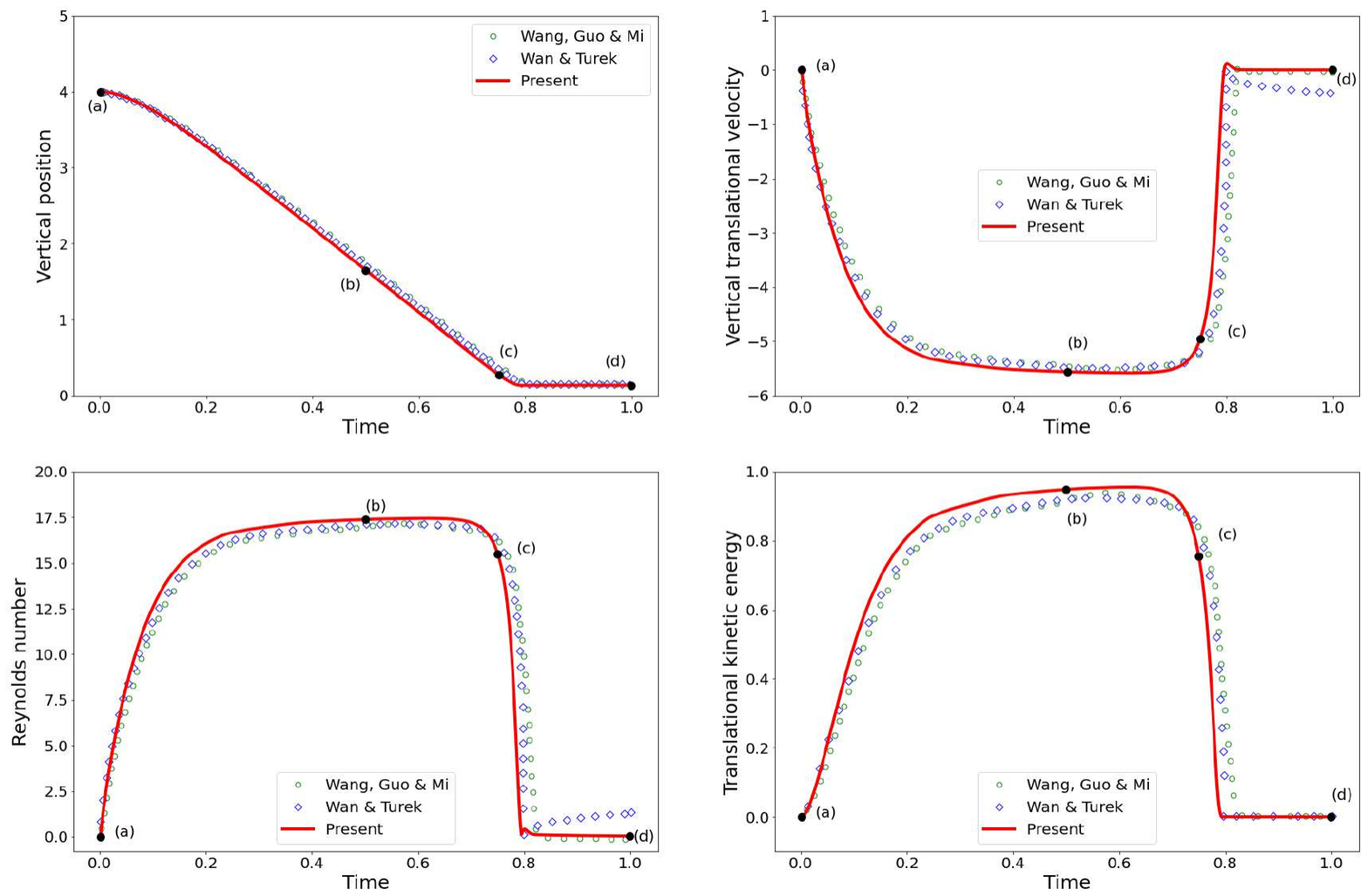}  
	
	\caption{Evolution of the vertical position of mass center, vertical translational 
	velocity, Reynolds number and translational kinetic energy. The red line are the present 
	results, the green dotted and blue dotted lines correspond to results respectively taken from 
	\cite{wang_drafting_2014} and \cite{wan_direct_2006}. The streamlines and body position at four 
	different time instants are represented.}

	\label{Fig:validation2d}
\end{figure}

For the validation in $2$D, we consider a rigid body of radius $r = 0.125$cm and density 
$\Density_{\Solid} = 1.25 \frac{\mbox{g}}{\mbox{cm}^3}$ 
is initially located at $(1$cm$,4$cm$)$ in a channel of width $2$cm and height $6$cm.  
Under the effect of gravity $g = 981 \frac{\mbox{cm}}{\mbox{s}^2}$ the particle is 
falling towards the bottom of a channel filled with a fluid of density $\rho_f = 1 \frac{\mbox{g}}{\mbox{cm}^3}$ 
and viscosity $\mu = 0.1 \frac{\mbox{g}}{\mbox{cm s}}$.
Both fluid and particle are at rest at time $t=0.0$s. The mesh size and time step used 
in this simulation are respectively fixed at $h = 0.01$cm and $\Delta t = 0.001$s. 
The width of the collision zone is set to $\rho = 0.015$cm. One can use a stiffness 
parameter for particle-fluid domain interaction $\epsilon_{\Fluid} < \frac{h^2}{2}$, for the results 
shown in the figures we used $\epsilon_{\Fluid} = 5*10^{-6}$. 
In figure~\ref{Fig:validation2d}, we plot the time evolution of four quantities 
and compare them to literature results Wan and Turek \cite{wan_direct_2006}, Wang, 
Guo and Mi \cite{wang_drafting_2014}. 
These quantities are the vertical coordinate of the center of mass  $y^{cm}$, 
the vertical translational velocity $v_y$, the Reynolds number $Re$ and the translational 
kinetic energy $E_t$ defined by:
$$
Re = \frac{2r \Density_{\Solid}\sqrt{v_x^2 + v_y^2}}{\mu} \quad \mbox{and} \quad E_t = 0.5\pi r^2 \Density_{\Solid} (v_x^2+v_y^2),
$$
where $v_x$ the horizontal translational velocity.
One can observe that all results are in good agreement. The small differences 
before and after collision can be explained by different definitions of collision parameters and 
numerical methods.   

\begin{figure}
	
	\centering  

	\includegraphics[scale=0.25]{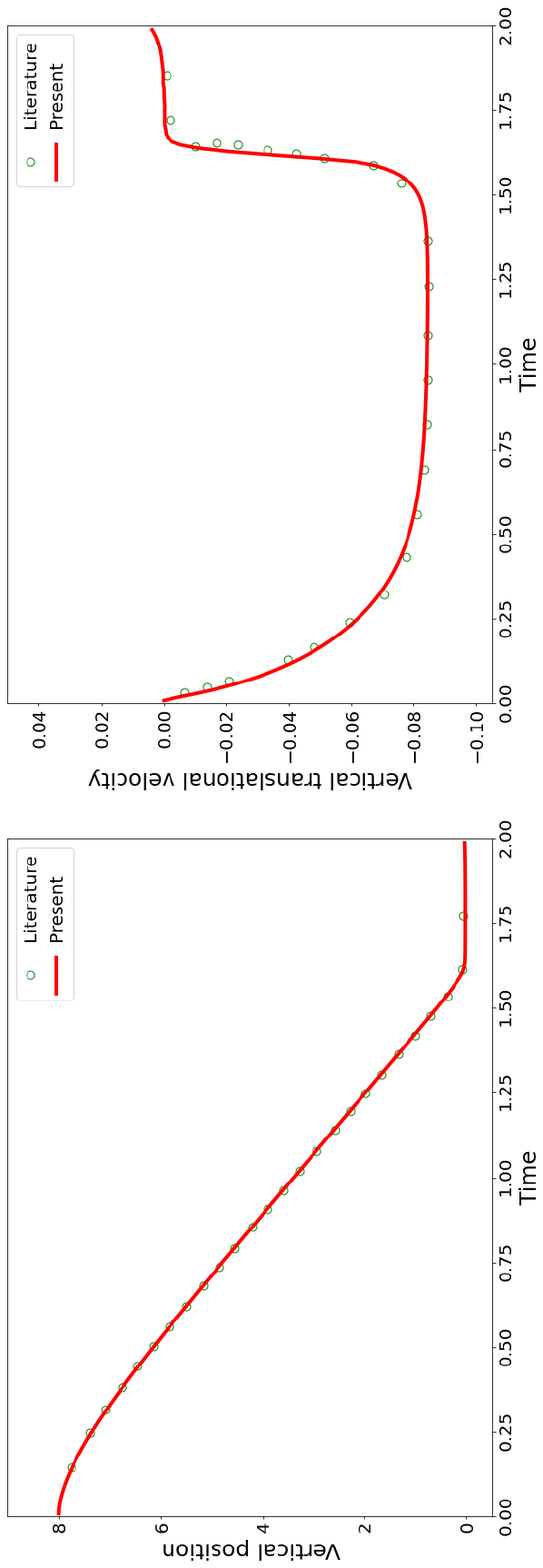}
	\caption{Comparison of the vertical position of the center of mass and
	vertical translational velocity observed in present results, red line, to 
	literature results \cite{zamora_numerical_2019}, dotted line.}

	\label{Fig:validation3d}
\end{figure}

We perform the same simulation for the $3$D validation. Initially, a sphere of 
radius $r=0.75$cm and density $\rho_s = 1.12 \frac{\mbox{g}}{\mbox{cm}^3}$ is located at 
$(0 \mbox{cm},12.75\mbox{cm},0 \mbox{cm})$ in a computational box of dimensions 
$[-4\mbox{cm},4\mbox{cm}] \times [0\mbox{cm},15\mbox{cm}] \times [-4\mbox{cm},4\mbox{cm}]$. 
The box is filled with a fluid of density $\rho_f = 0.962 \frac{\mbox{g}}{\mbox{cm}^3}$ and 
viscosity $\mu = 1.13 \frac{\mbox{g}}{\mbox{cm s}}$. The mesh size and the time 
step are set to $h = 0.1$cm and $\Delta t = 0.01$s. Regarding the lubrication force parameters, 
we fixed the width of the collision zone to $\rho = 0.2625$cm and the stiffness 
parameter to $\epsilon_{\Fluid} = 5*10^{-6}$. The comparison of the present results to 
literature results \cite{zamora_numerical_2019} is given in figure \ref{Fig:validation3d}. We obtain 
the same results for the vertical coordinate of the center of mass and the vertical 
translational velocity.  \newline 

\section{Applications} \label{applications}

\subsection{Isolated object}

In this subsection we simulate the interaction between a complex shaped body falling under 
the effect of gravity and the boundary of the computational domain filled with an incompressible 
Newtonian fluid in two and three dimensions.  \newline 

\begin{figure}[h!]
    \begin{center}
	\begin{subfigure}{0.49\textwidth}
		\centering
		\includegraphics[scale=0.35]{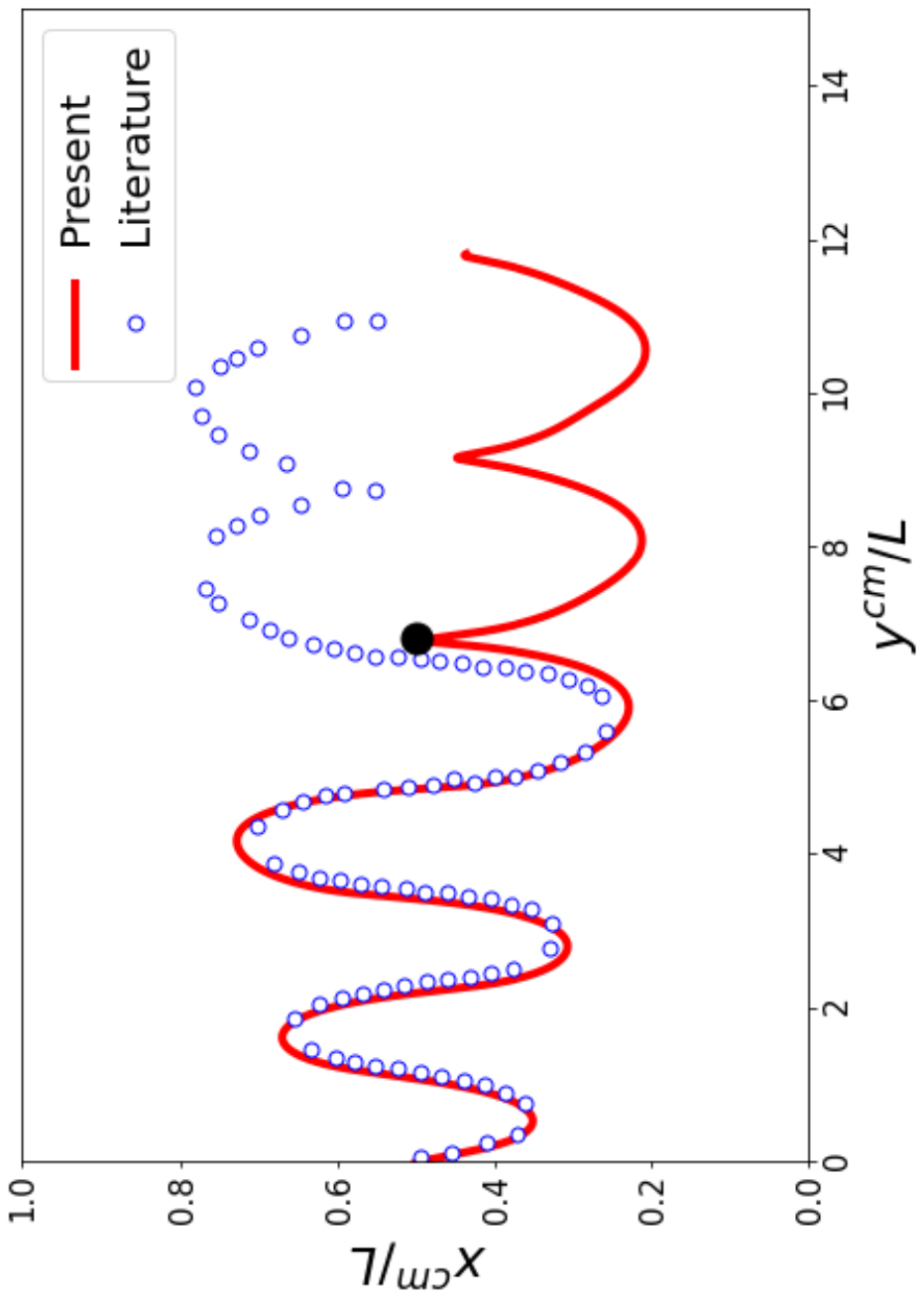}
		\caption{Trajectory}
		\label{Fig:ellipse_1}
	\end{subfigure}
	\begin{subfigure}{0.49\textwidth}
		\centering
		\includegraphics[scale=0.35]{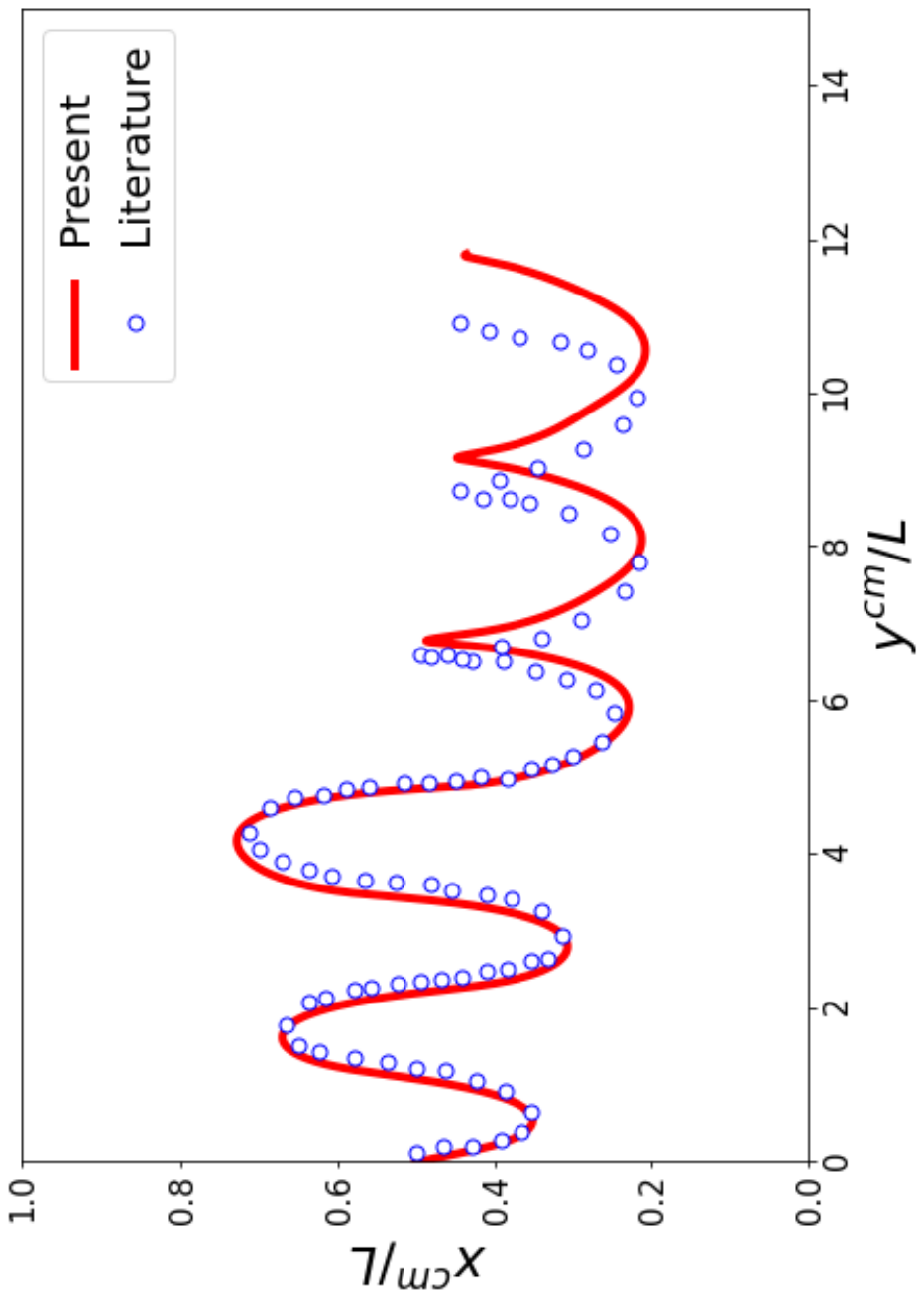}
		\caption{Trajectory with symmetric results}
		\label{Fig:ellipse_2}
	\end{subfigure}
	\end{center}

	\caption{The left graph plots the trajectory of the ellipse center of mass. 
	The red line corresponds to the present results and the blue dotted line represents 
	the literature results \cite{xia_flow_2009}. The right graph shows the same trajectory but 
	compares it to rotations happening on the other wall.}

	\label{Fig:ellipse}
\end{figure}

For the first case, we consider a two-dimensional ellipse located in a channel of width $L = \frac{16}{130}$cm 
and infinite length. The channel contains a fluid of density 
$\Density_{f} = 1.195\frac{\mbox{g}}{\mbox{cm}^3}$ and viscosity 
$\mu = 0.305\frac{\mbox{g}}{\mbox{cm s}}$. The ellipse's long and short axis are 
respectively set to $a = 0.1$cm and $b = 0.05$cm. Its density is equal to 
$\Density_{\Solid} = 1.35\frac{\mbox{g}}{\mbox{cm}^3}$. Due to its initial orientation, 
set to $\theta = \frac{\pi}{3}$ with respect to the horizontal axis, the ellipse 
will collide with the right and left walls. 
To model these collisions the length of the lubrication zone is fixed at $\rho = 0.05$cm 
and a stiffness parameter of $\epsilon_{\Fluid} = 10^{-5}$ is used. The simulation is performed 
using a mesh size $h = 0.001$cm and a time step $\Delta t = 0.001$s. 
Figure \ref{Fig:ellipse_1} shows the trajectory of the ellipse comparing it to literature \cite{xia_flow_2009}.
This trajectory allows to observe that the ellipse first performs oscillations between the 
right and left wall before reaching a horizontal position, at the instant given by the black point. 
Then, it performs rotations close to a single wall. In present results, these rotations 
are taking place near the left wall, but in the reference article, they happen near the right wall. 
Given that the ellipse is in a horizontal position at the beginning of these rotations, 
any small difference in the collision model or the numerical resolution techniques can 
explain this difference.
To show that the present results remain close to literature, we consider in 
figure \ref{Fig:ellipse_2} the symmetric trajectory for the rotations of \cite{xia_flow_2009}. 
This test case allows us to validate our implementation for the case of arbitrarily shaped bodies.\newline

\begin{figure}[h!]
	\centering
	\includegraphics[scale=0.35]{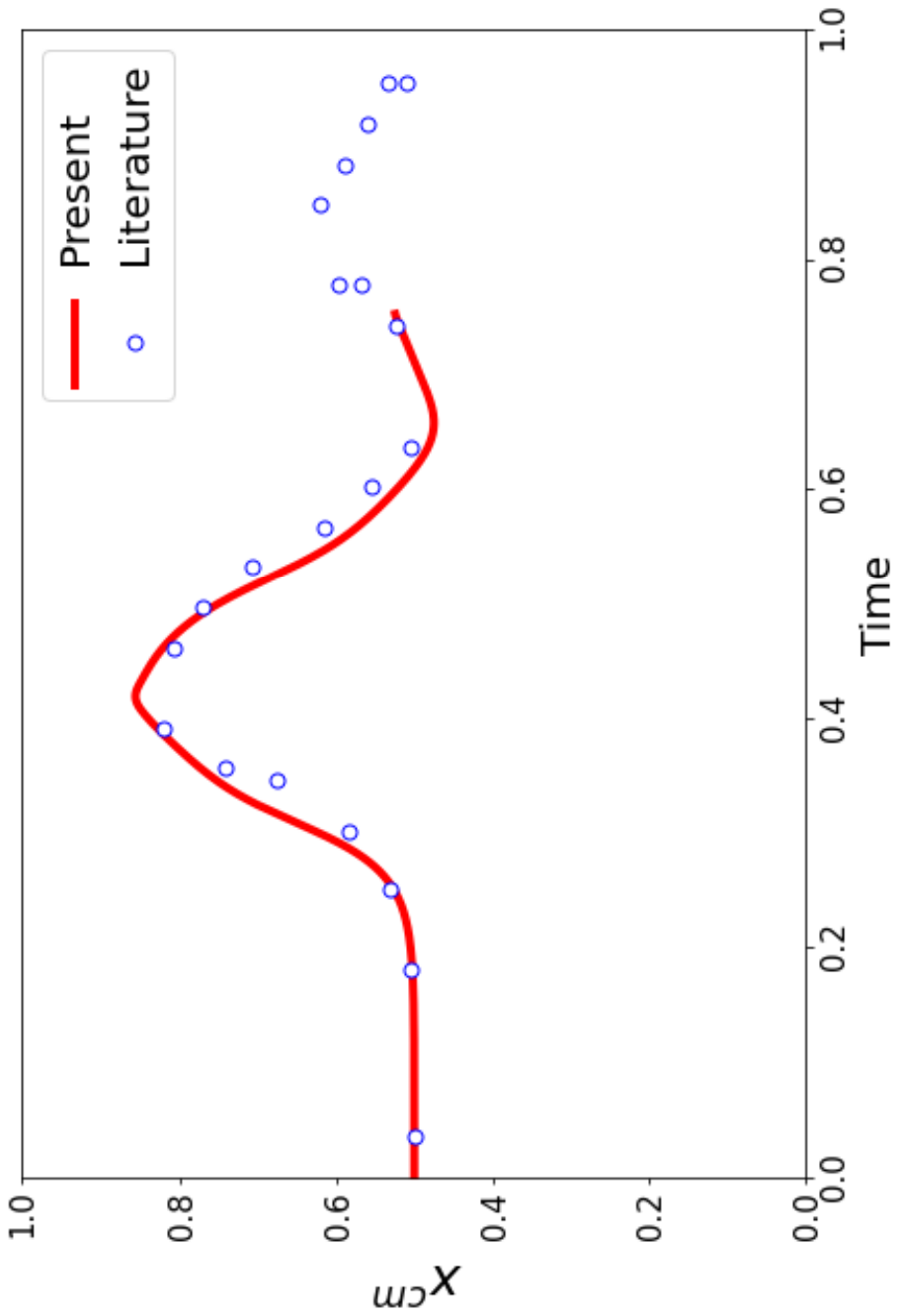}

	\caption{Evolution of the horizontal position of the center of the ellipsoid. The 
	dotted line shows the literature results \cite{pan_direct_2002} and the red line 
	represents the present results.}

	\label{Fig:ellipsoid}
\end{figure}

The second case represents the simulation of a three-dimensional ellipsoid of 
density $\Density_{\Solid} = 1.25\frac{\mbox{g}}{\mbox{cm}^3}$ and axis 
of length $a = 0.4$cm and $b = c = 0.2$cm falling in a rectangular 
computational domain filled with a fluid of density 
$\Density_{f} = 1\frac{\mbox{g}}{\mbox{cm}^3}$ and viscosity 
$\mu = 0.01 \frac{\mbox{g}}{\mbox{cm s}}$. The dimensions of the domain are set to 
$[0\mbox{cm},1\mbox{cm}]\times [0\mbox{cm},8\mbox{cm}]\times [0\mbox{cm},0.4\mbox{cm}]$
and at $t=0$s the ellipsoid is located at $(0.5\mbox{cm}, 6 \mbox{cm}, 0.2\mbox{cm})$ 
with its long axis oriented in vertical direction. The mesh size and time step of this 
simulation are $h=0.0125$cm and $\Delta t = 0.001$s. The interactions between the ellipsoid 
and the boundary of the computational domain are modelled by defining a collision zone 
of width $\rho = 0.04$cm and a stiffness parameter $\epsilon_{\Fluid} = 3*10^{-6}$. The 
evolution of the horizontal position of the sphere is shown in figure \ref{Fig:ellipsoid}, 
and is compared to literature results \cite{pan_direct_2002}. The trajectories 
before and after the first interaction between the ellipsoid and the boundary 
of the computational domain, at time $t \approx 0.45$s, are in good agreement.  

\subsection{Multiple objects}

\paragraph{Two disks falling in an incompressible fluid}

\begin{figure}[h!]
    \begin{center}
		\begin{subfigure}{0.22\textwidth}
			\centering
			\includegraphics[scale=0.3]{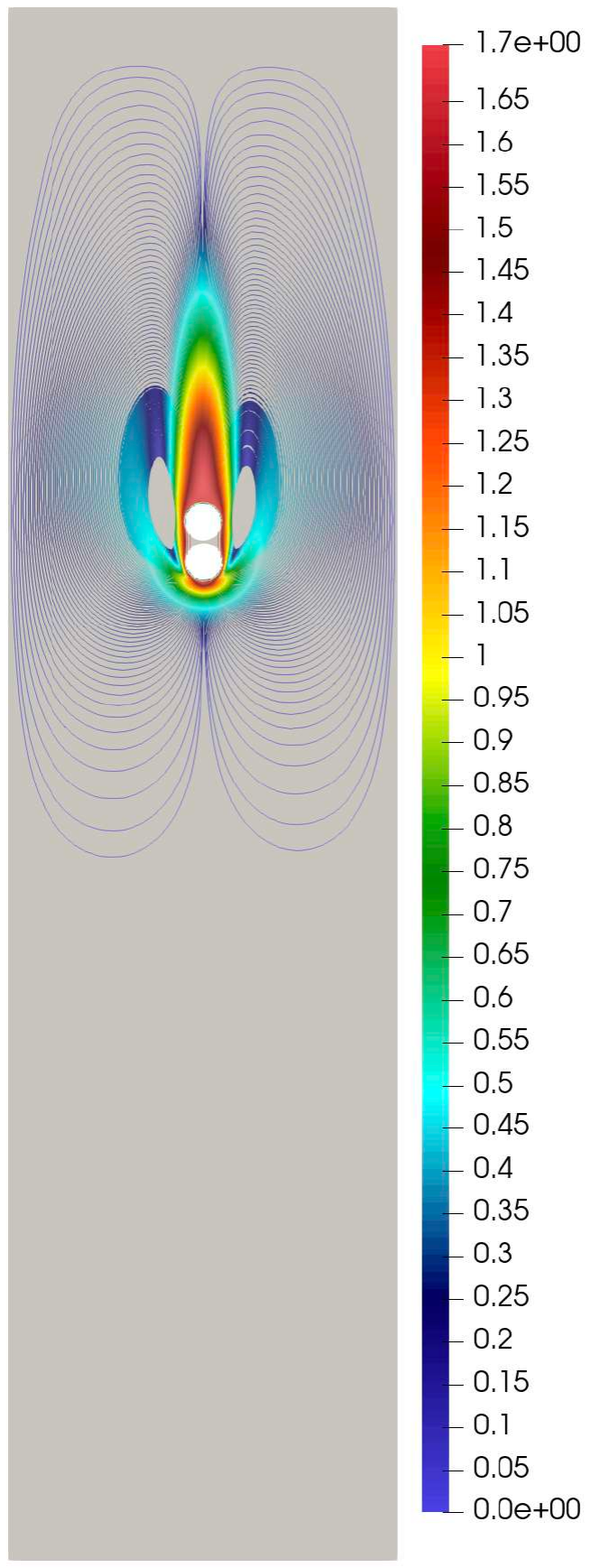}
			\caption{$t=1.5 s$}
		\end{subfigure}
		\begin{subfigure}{0.22\textwidth}
			\centering
			\includegraphics[scale=0.3]{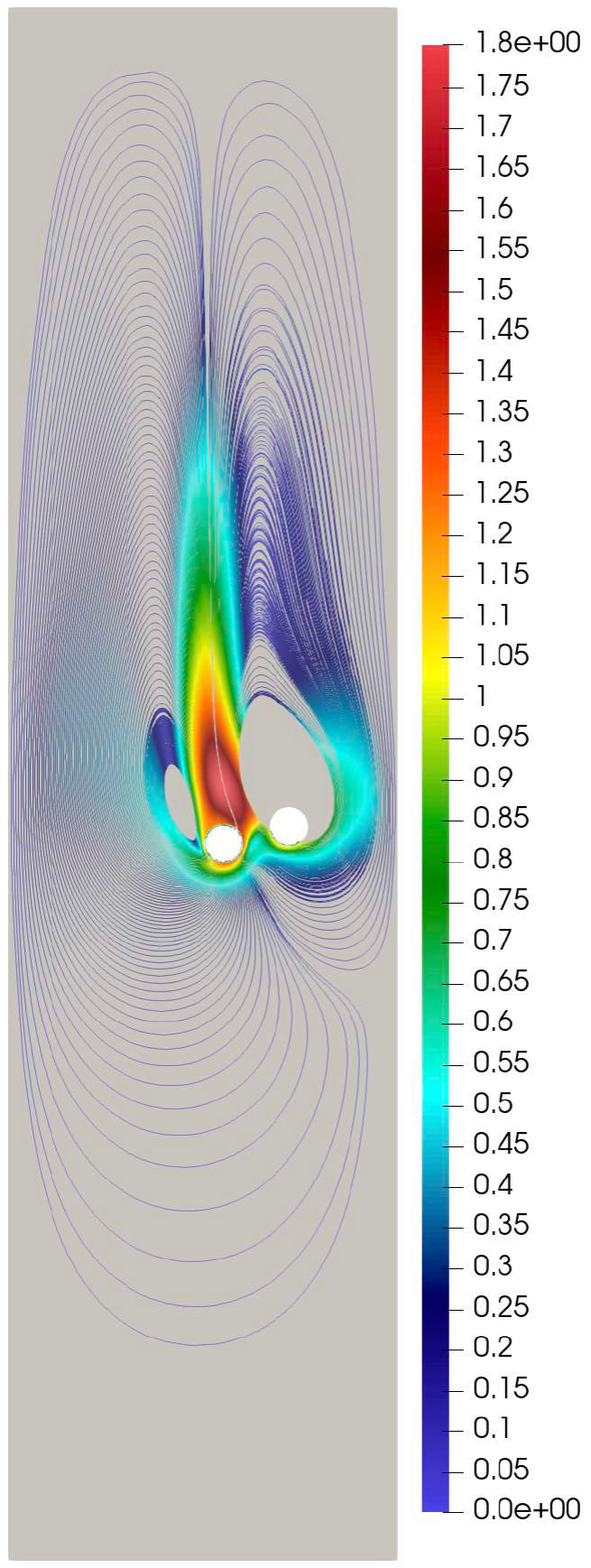}
			\caption{$t=2.5 s$}
		\end{subfigure}
		\begin{subfigure}{0.22\textwidth}
			\centering
			\includegraphics[scale=0.3]{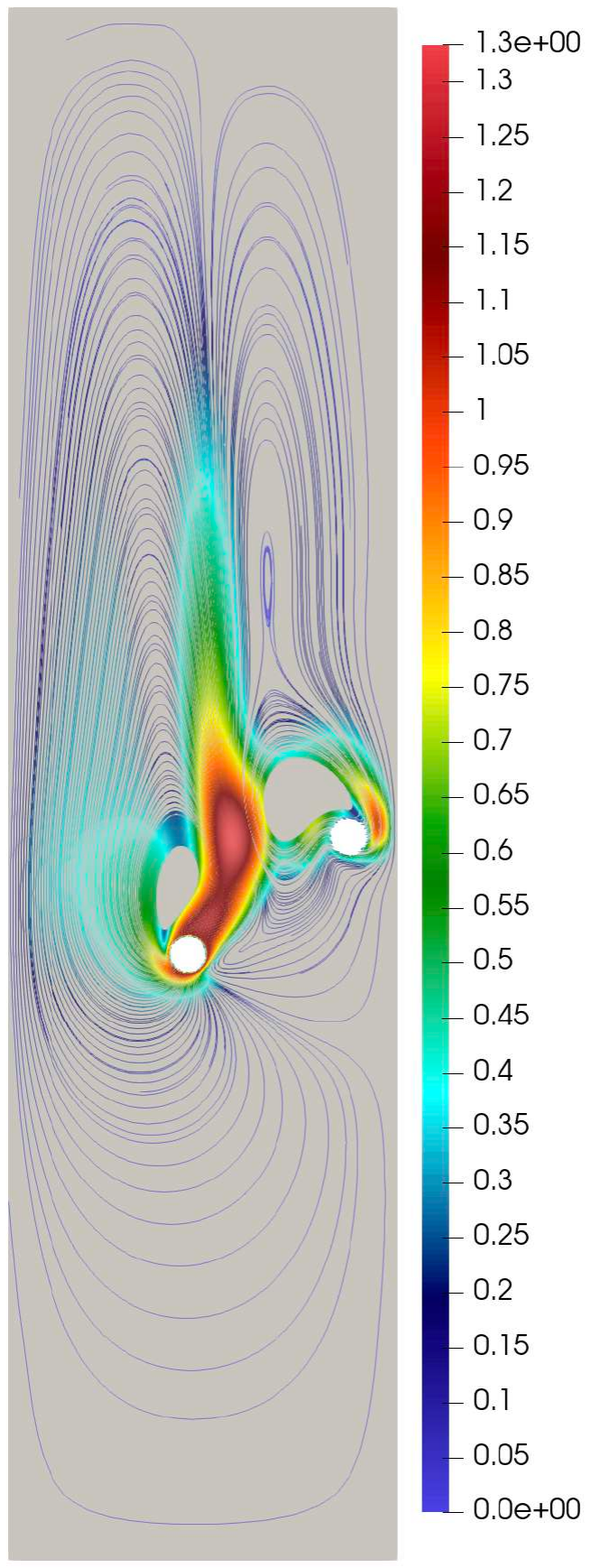}
			\caption{$t=3.0 s$}
		\end{subfigure}
		\begin{subfigure}{0.24\textwidth}
			\centering
			\includegraphics[scale=0.3]{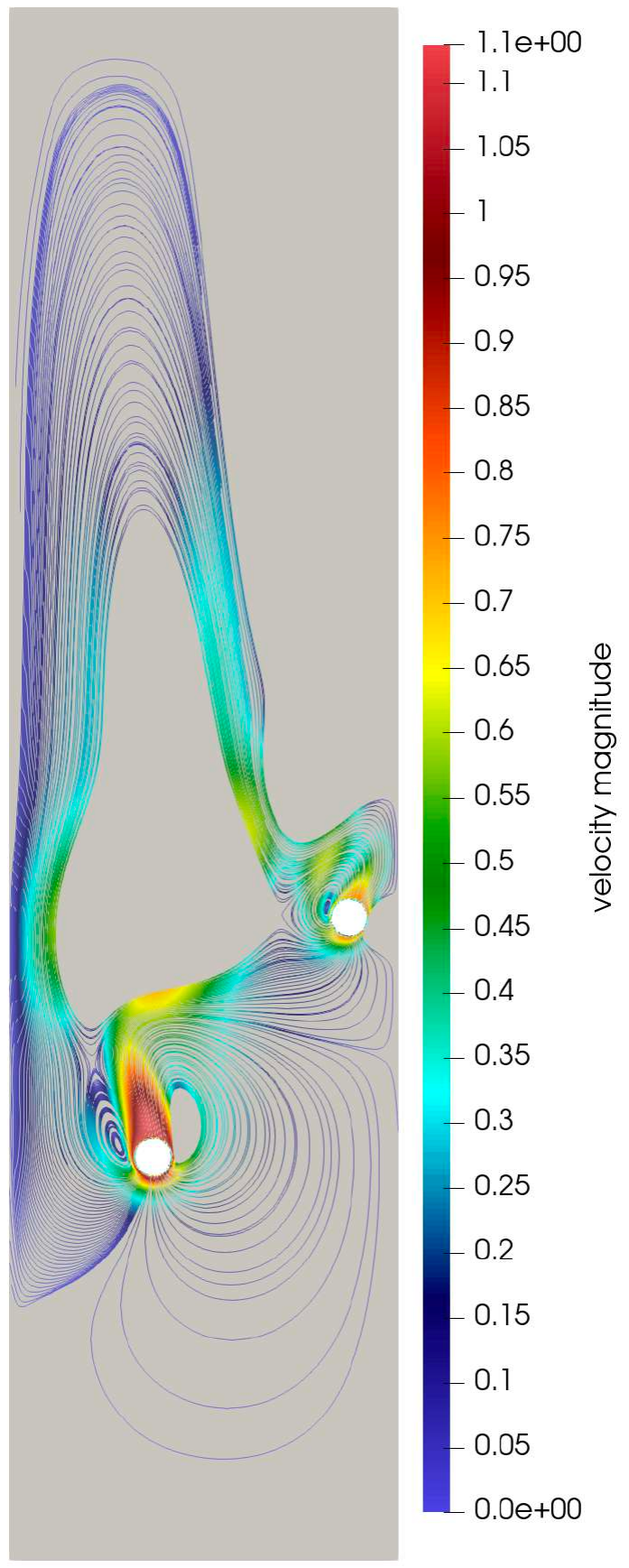}
			\caption{$t=4.0 s$}
		\end{subfigure}
		\end{center}

		\centering  
		\includegraphics[scale=0.17]{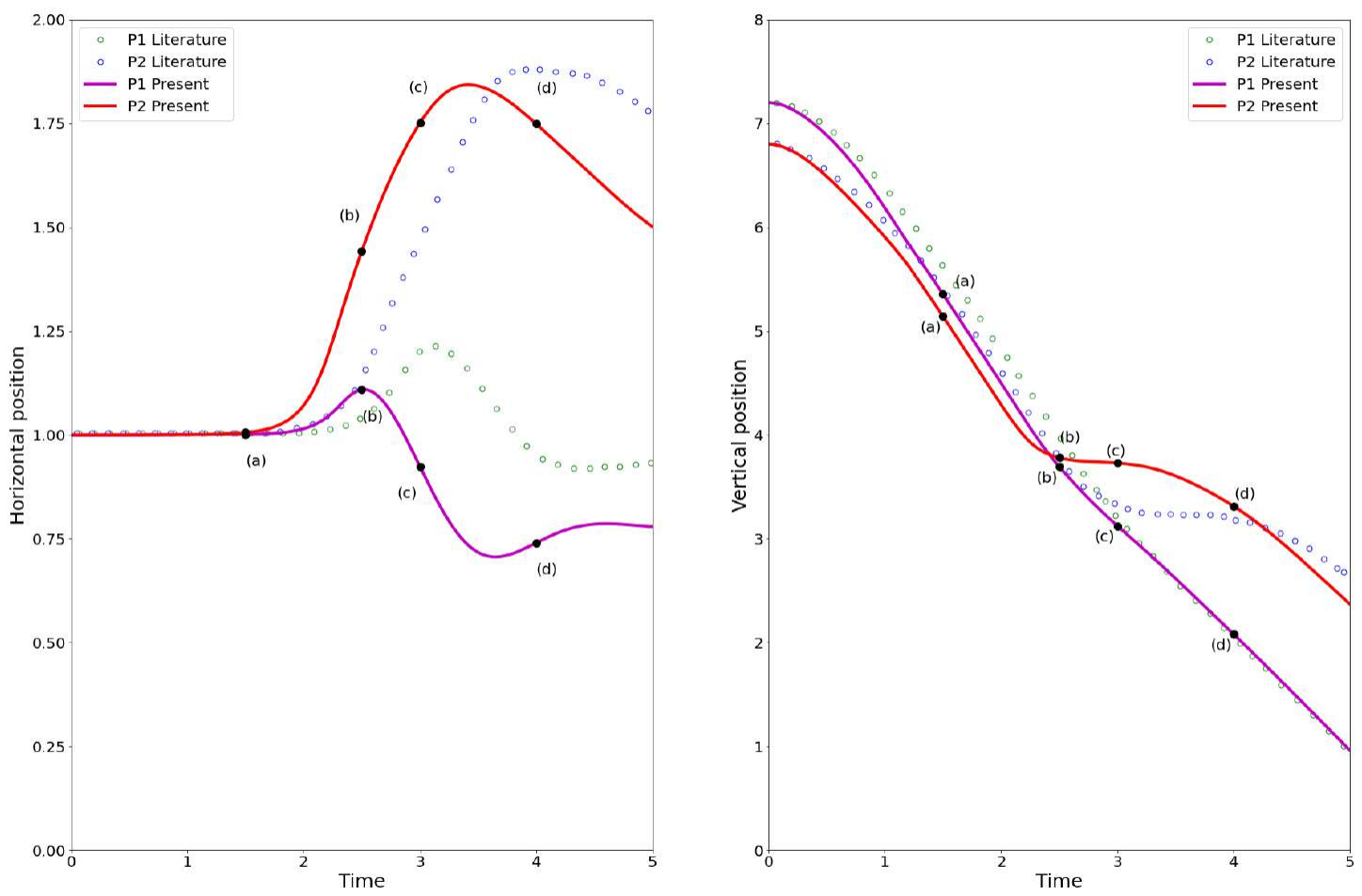}

	\caption{The upper figures show the position of the two disks as well as the 
	streamlines at four different time instants. 
	The lower graphs compare the evolution of the disks horizontal and vertical 
	position to literature \cite{feng_immersed_2004}. The solid lines correspond 
	to present results and the dotted one to literature results.}
	\label{Fig:twoDisks}
\end{figure}

Collisions between two spherical bodies placed on a same vertical line and falling 
under the effect of gravity, in a channel filled with incompressible Newtonian fluid, 
are known as the \textit{drafting, kissing and tumbling} phenomenon. 
We consider two disks of radius $r = 0.1$cm and density $\Density_{\Solid} = 1.01 \frac{\mbox{g}}{\mbox{cm}^3}$ 
placed in a fluid domain of dimensions $[0\mbox{cm},2\mbox{cm}] \times [0\mbox{cm},8\mbox{cm}]$. 
To have a slight asymmetry, necessary to speed up the insurgence of the phenomenon, the coordinates 
of the upper disk  $\Solid 1$ center are fixed at $(0.999 \mbox{cm},7.2 \mbox{cm})$ and those of the 
lower disk $\Solid 2$ center at $(1.0 \mbox{cm},6.8 \mbox{cm})$. The fluid density and viscosity 
are respectively given by $\rho_f = 1.0 \frac{\mbox{g}}{\mbox{cm}^3}$ and $\mu = 0.01 \frac{\mbox{g}}{\mbox{cm s}}$. 
The simulation is performed on the time interval $[0\mbox{s},5\mbox{s}]$ with a time step equal 
to $\Delta t = 0.001$s and a mesh size set to $h = 0.01$.
To model the \textit{drafting, kissing and tumbling} phenomenon as in \cite{feng_immersed_2004}, one fixes the width of 
the collision zone at $\rho = 0.0225$cm and uses stiffness parameters for the disk-disk and 
disk-fluid domain collisions given by $\epsilon = 7*10^{-5}$ and $\epsilon_{\Fluid} = 5*10^{-5}$. Figure \ref{Fig:twoDisks} 
shows the results of the simulation. When both 
disks are next to each other, the motion of the lower disk reduces the fluid pressure 
behind it and therefore the 
resistance of the fluid for the upper disk. As a result, the upper body falls faster 
until it collides with the lower one. The collision forces cause the separation
of the disks which then move in opposite directions. Figure \ref{Fig:twoDisks} 
also compares the vertical and horizontal position of the two disks over time to 
literature \cite{feng_immersed_2004}. It can be observed that the results of the 
collision phase have same behavior.  In the present simulation, the bodies fall faster 
to the bottom, therefore the separation phase takes place sooner than in 
literature results. This explains the differences on the graphes.\newline

\paragraph{$100$ disks falling in an incompressible fluid}

\begin{figure}[h!]
	
	\begin{center}
	\begin{subfigure}{0.2\textwidth}
		\centering
		\includegraphics[scale=0.23]{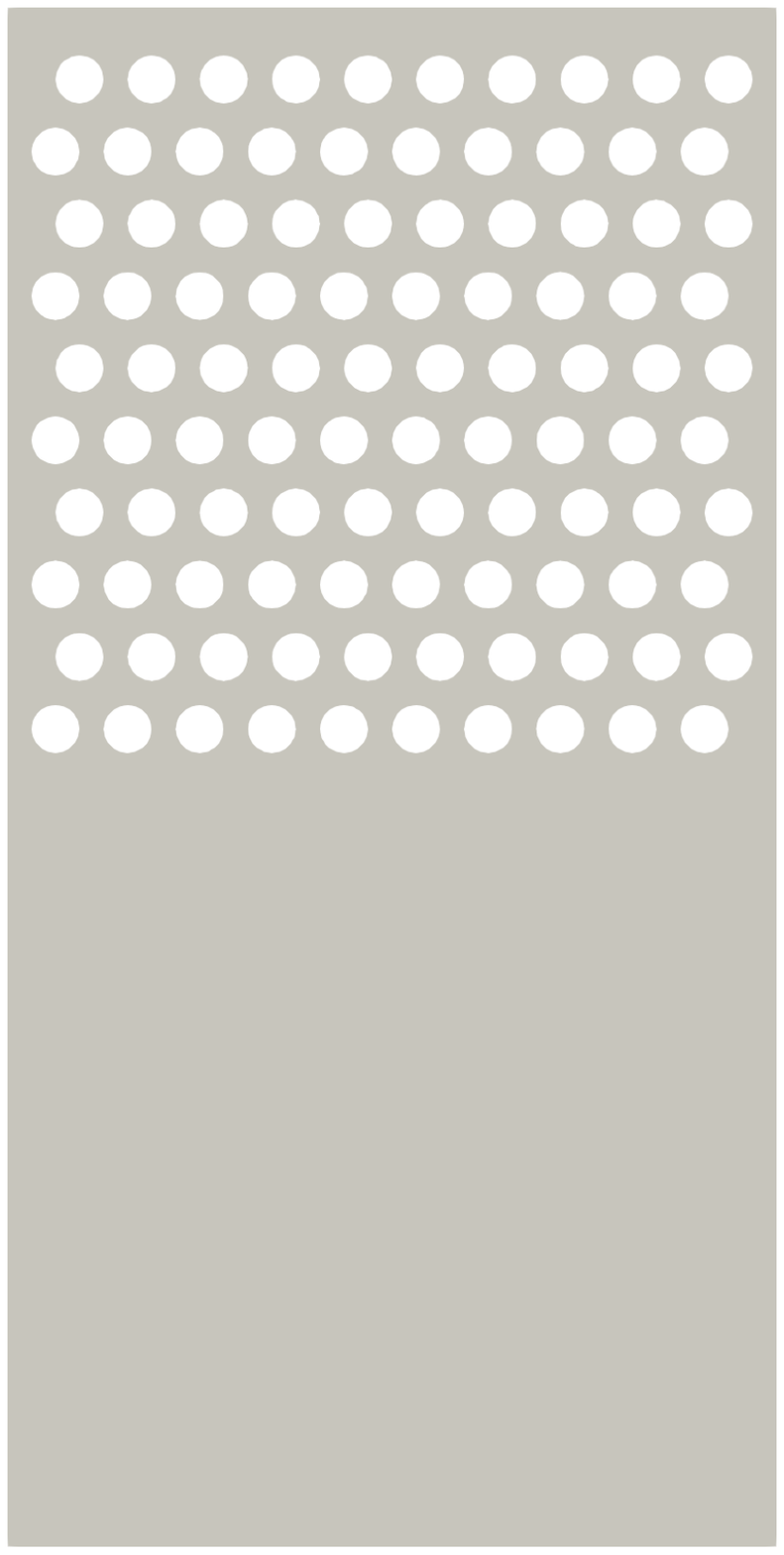}
		\caption{$t=0.0 s$}
	\end{subfigure}
	\begin{subfigure}{0.27\textwidth}
		\centering
		\includegraphics[scale=0.23]{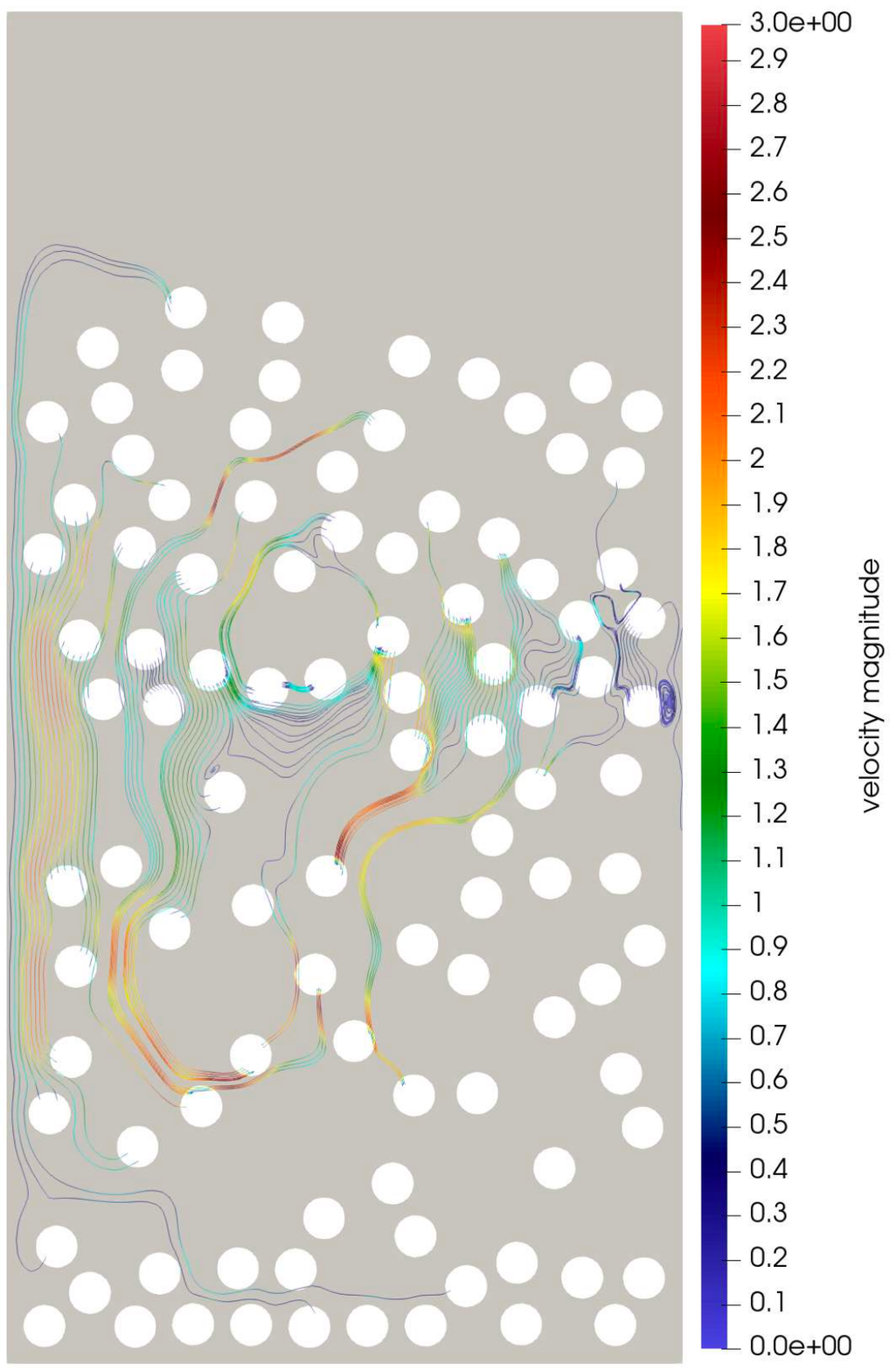}
		\caption{$t=1.5 s$}
	\end{subfigure}
	\begin{subfigure}{0.27\textwidth}
		\centering
		\includegraphics[scale=0.23]{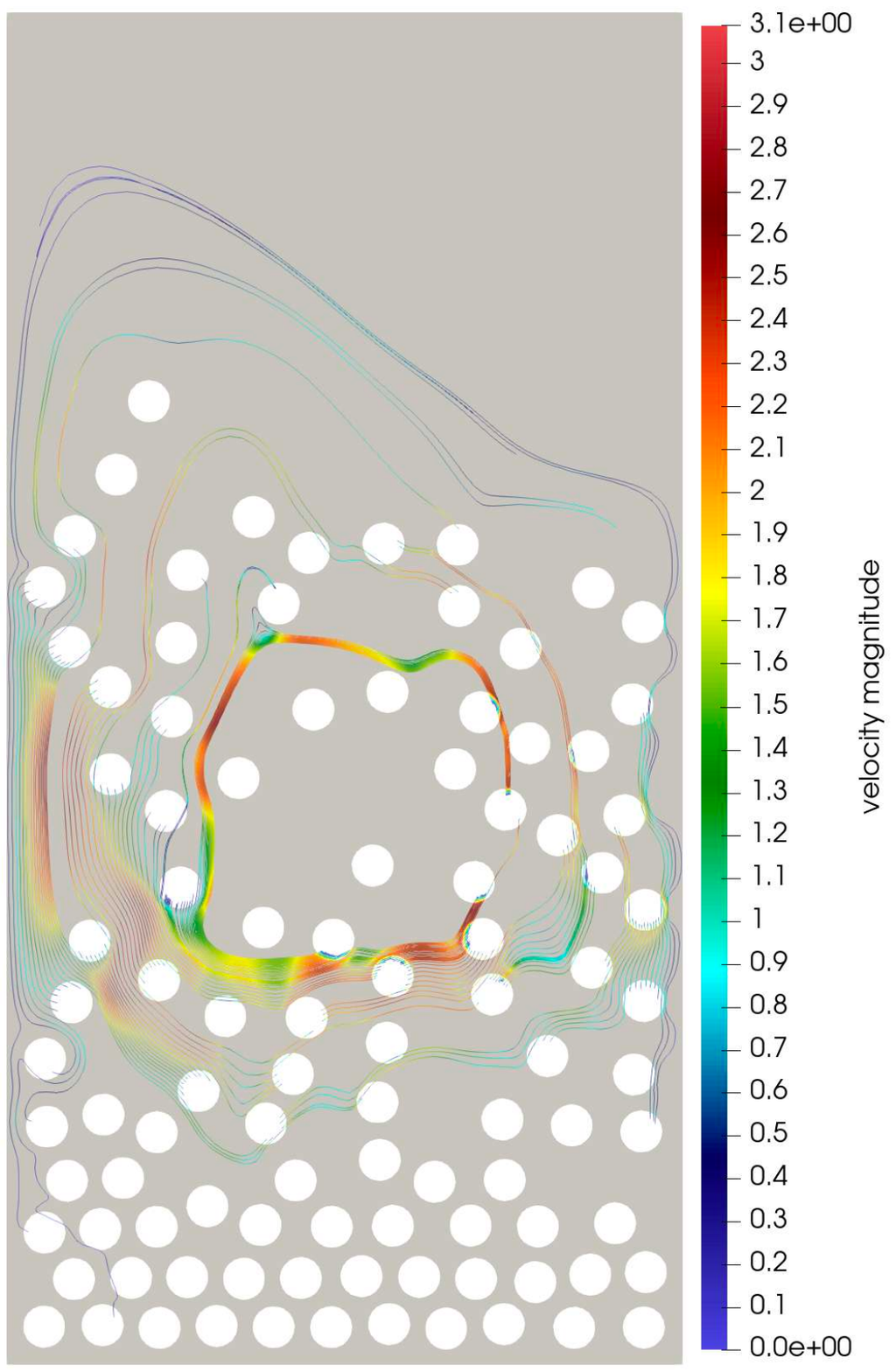}
		\caption{$t=2.0 s$}
	\end{subfigure}
	\begin{subfigure}{0.2\textwidth}
		\centering
		\includegraphics[scale=0.23]{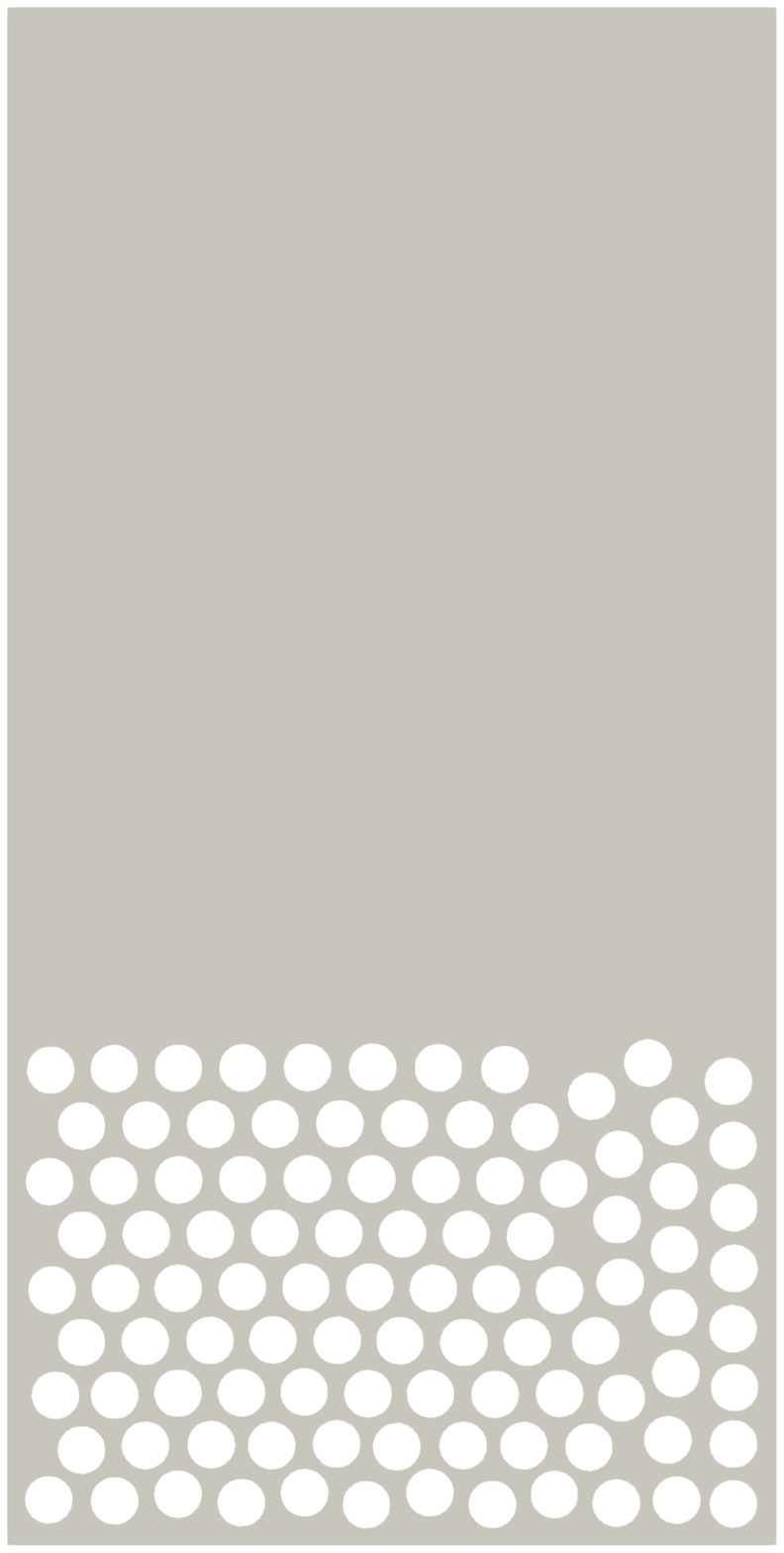}
		\caption{$t=3.5 s$}
	\end{subfigure}
	\end{center}

	\caption{Streamlines and position at four different time instants of 100 disks falling under the effect of gravity in an incompressible Newtonian fluid.}

	\label{Fig:100disks}
\end{figure}

We simulate the interactions between $100$ circular 
particles in two dimensions. The particles of radius $r = 0.03125$cm and density 
$\Density_{\Solid} = 1.1 \frac{\mbox{g}}{\mbox{cm}^3}$ are placed in a vertical channel of
 height $2$cm and width $1$cm which is filled with a fluid of density 
 $\rho_f = 1. 0 \frac{\mbox{g}}{\mbox{cm}^3}$ and viscosity $\mu = 0.01 \frac{\mbox{g}}{\mbox{cm s} }$. 
 At the beginning of the simulation, the fluid and the particles are at rest. Then the 
 particles move towards the bottom of the domain under the effect of gravity. 
 The initial configuration is given by the first figure of \ref{Fig:100disks}. 
 We run the simulation on a time interval $[0\mbox{s},5.5\mbox{s}]$ and use a time 
 step fixed at $\Delta t = 0.002$s. The mesh size is given by $h = 0.01$cm.
The collision force is defined on range $\rho = 0.025$cm. According to simulations 
of multiple objects interactions in the literature \cite{pan_fluidization_2002}, the force intensity must be 
important to prevent the objects from overlapping. For this reason, we use 
stiffness parameters fixed at $\epsilon = 5*10^{-6}$ and $\epsilon_{\Fluid} = 10^{-7}$. 
Figure \ref{Fig:100disks} shows the position of the $100$ particles at four time instants. 
Due to the collision forces, the particles near the domain walls fall slower towards 
the bottom than the central particles. During the simulation the particles 
settle one onto the other on the bottom of the channel. At final time, almost all particles 
are in a stationary position and packed in a hexagonal lattice. 
Article \cite{juarez_numerical_nodate} illustrates the same simulation but with different values of the parameters. 

\paragraph*{Two objects in a flow}

For this application we consider the motion of two circular particles in a symmetric 
stenotic artery in two dimensions. In contrast with the previous test cases, the particles 
do not fall under the effect of gravity, but the motion of the particles 
is due to a pressure difference between the inlet and the outlet. This simulation 
is presented in \cite{shu_particulate_2010} and \cite{li_lattice_2004}.\newline 

\begin{figure}[h!]
    \centering
    \includegraphics[scale=0.35]{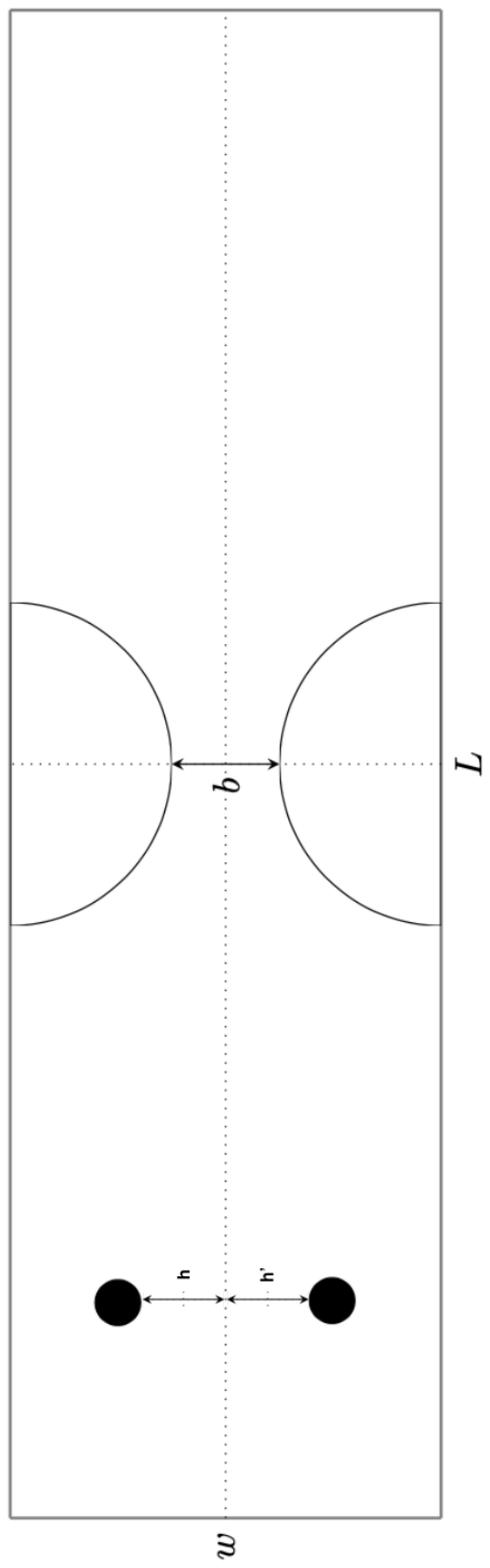}
    \caption{Geometry of symmetric stenotic artery}
    \label{fig:artery}
\end{figure}

\begin{figure}[h!]
    \centering
    \includegraphics[scale=0.55]{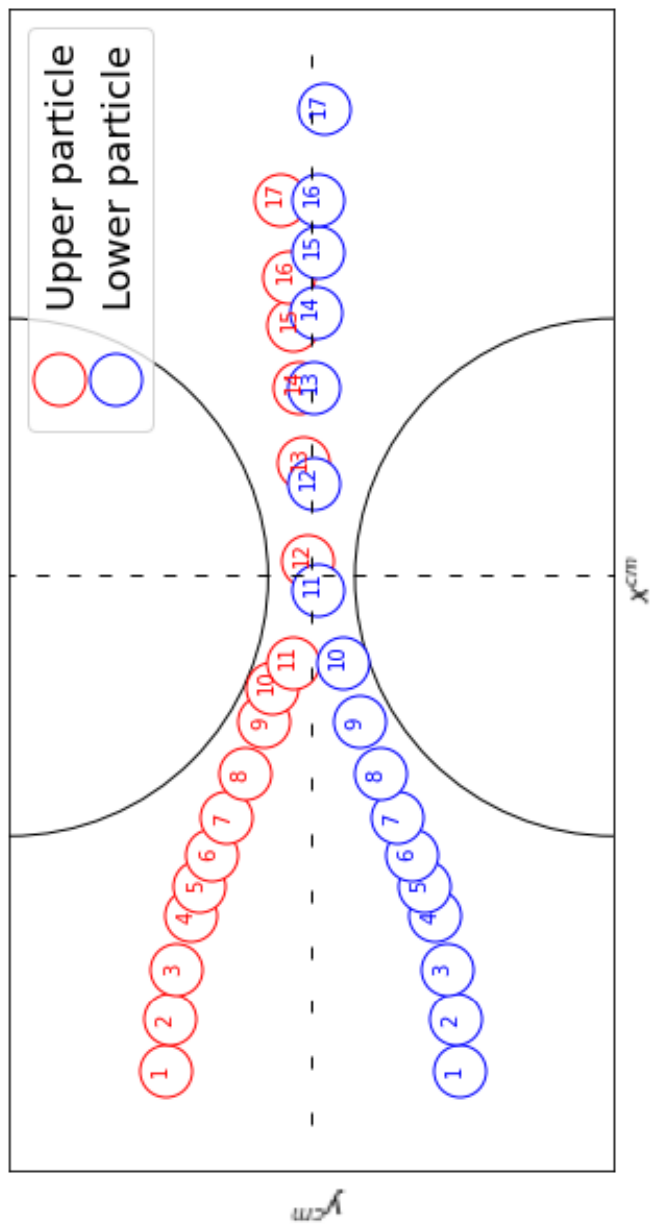}
    \caption{Position of two particles at different time instants.}
    \label{fig:artery_pos}
\end{figure}

The geometry, given in figure \ref{fig:artery}, is a channel of length 
$L = 32d$ and width $w = 8d$. 
The centerlines are represented by dotted lines. 
The diameter $d$ of the particles is set to $d = 8.5 * 10^{-4}$cm. To create the 
stenosis, the authors of \cite{shu_particulate_2010} and \cite{li_lattice_2004} add 
two symmetrical protuberances shaped like semicircles. The radius of these protuberances 
depends on the width of the stenosis throat, which is fixed at $b = 1.75d$. The fluid viscosity is 
set to $\mu = 0.01\frac{\mbox{cm}^2}{\mbox{s}}$ and both fluid and particle 
densities are equal to $\rho_f = \Density_{\Solid}= 1\frac{\mbox{g}}{\mbox{cm}^3}$.
At time $t=0$, the fluid and the particles are at rest. Then the particles move towards 
the protuberances due to the pressure difference $\Delta p = 541 \mbox{Pa}$ between 
the inlet and outlet of the channel. The two 
particles are initially located at $8d$ to the left of the vertical centerline and asymmetrically 
with respect to the horizontal centerline: $h = 2d + \frac{d}{4000}$ and $h' = 2d$. This slight 
asymmetry allows them to cross the stenosis throat. The mesh size and time step of 
the simulation are respectively set to $h = 6.07 * 10^{-5}$cm and 
$\Delta t = 10^{-5}$s. The time interval is $[0.0\mbox{s},0.01\mbox{s}]$.
Collision forces are applied when the particles are close to the protuberances. 
The width of the collision zone is fixed at $\rho = 0.0003$cm and the stiffness 
parameters of the particle-particle and particle-domain interactions are given by 
$\epsilon = 10^{-11}$ and $\epsilon_{\Fluid} = 7.5*10^{-12}$. 
Figure \ref{fig:artery_pos} shows the trajectory of the two particles: 
they have a similar motion (snapshots $1$-$9$) until they get close 
to the stenosis throat (snapshot $10$). Since the width of the stenosis throat 
is too small to allow both particles to pass simultaneously, the upper particle 
stops and changes direction due to collision forces. Snapshot $11$ shows that it 
moves back to allow the lower particle passing the throat. Finally, the upper 
particle follows the lower particle through the artery (snapshots $11$-$17$). 
The same particle behavior is present in the cited articles. The particle trajectory 
after passing the stenosis throat depends on the initial configuration and 
the definition of collision parameters.

\begin{figure}[h!]
    \centering
    \includegraphics[scale=0.59]{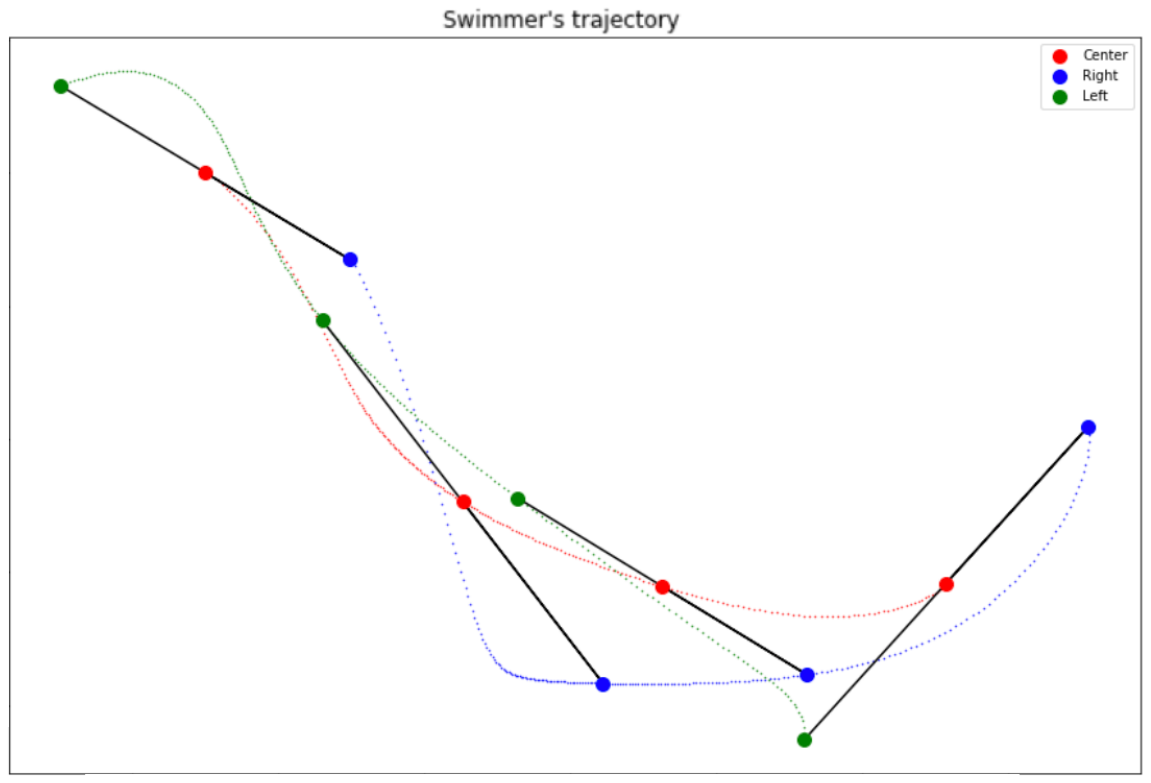}
    \caption{Behavior of the three-sphere swimmer close to the boundary of 
	the computational domain.}
    \label{fig:3ss}
\end{figure}

\paragraph*{Three-sphere swimmer in a Stokes flow}
For the last test case, we consider the tilted three-sphere swimmer placed close to 
the boundary of a horizontal channel in a Stokes regime. The radius of the swimmer 
spheres is set to $r = 1$cm and the length of its rods to $l=10$cm. The height and 
width of the computational domain are fixed at $h = 150$cm and $w = 40$cm. The 
mass center of the central sphere is initially located at $(15 \mbox{cm}, 10 \mbox{cm})$, 
and the density of 
the swimmer is given by $\rho_{\Solid} = 0.1 \frac{\mbox{g}}{\mbox{cm}^3}$. 
The fluid has a density equal to $\rho_f = 1 \frac{\mbox{g}}{\mbox{cm}^3}$ and a viscosity 
of $\mu = 1 \frac{\mbox{cm}^2}{\mbox{s}}$. The mesh size is set to $h = 0.3$cm and 
the time step to $\Delta t = 0.1$s. Due to the orientation of the swimmer at 
$t=0$s and due to its swimming strategy, leading to a straight motion, the swimmer 
gets closer to the boundary of the channel. Once its right sphere 
arrives in the collision zone of width $\rho = 0.03$cm, lubrication forces, defined 
by a stiffness parameter set to $\epsilon_{\Fluid} = 5*10^{-8}$, are applied onto 
the swimmer, who starts to change direction. The swimmer rotates until its left sphere 
enters the collision zone, where the repulsive force applied on this sphere forces the swimmer 
to move upwards, away from the boundary.
This behavior of the three-sphere swimmer close to the boundary is shown in figure 
\ref{fig:3ss}. 

\section*{Reproducibility}

The validation benchmarks and all applications illustrated in this paper 
are available in a public GitHub repository \cite{feelpp}. All the results can be 
reproduced. 

\section{Conclusions} \label{conclusion}

In this paper, we presented a repulsive collision-avoidance model for arbitrary-shaped 
bodies, where the definitions of force direction and magnitude are based on the 
computation of the distance function using the narrow-band fast marching method. 
The model is validated against literature results, with which we find good agreement. 
Various applications of two and three-dimensional cases are illustrated.  
In the appendix, we present a tool for the insertion of solid bodies into meshed domains, which 
allows simulations where the conforming treatment of the fluid-solid interface is required. 
Determining the optimal values for the collision parameters is still an open problem 
which could be tackled by anticipating the moments of collision and controlling the 
displacement of the bodies. The extension of this work to elastic bodies and collisions with 
multiple contact points will be the subject of an upcoming paper. 

\appendix

\section{Insertion of arbitrary solids into discretized fluid domains}

In biological processes, such as particle transport in blood vessels, solid bodies 
move in geometrically complex domains. The reconstruction of such environments 
for numerical simulations is often image-based rather than relying on a CAD 
description, which complicates the mesh-fitted definition of solid bodies inside 
these environments. Hence, we have developed a tool allowing the insertion of 
arbitrary bodies in meshed domains with the identification of their geometric 
and material properties.
The tool relies on the MMG library \cite{mmg} feature that allows remeshing from a 
level-set function: given the initial triangulated environment and a level-set 
description of the solid body, it outputs a new triangulation that includes 
the solid, and where the interface between the two is conforming. 
In addition to this, we are able to insert multiple objects simultaneously, 
with prescribed position and orientation, and ensure that physical markers 
associated to the fluid and the solids are correctly transferred onto the 
new mesh. \newline 

\begin{figure}[h!]
    \centering
	\centerline{\includegraphics[scale=0.35]{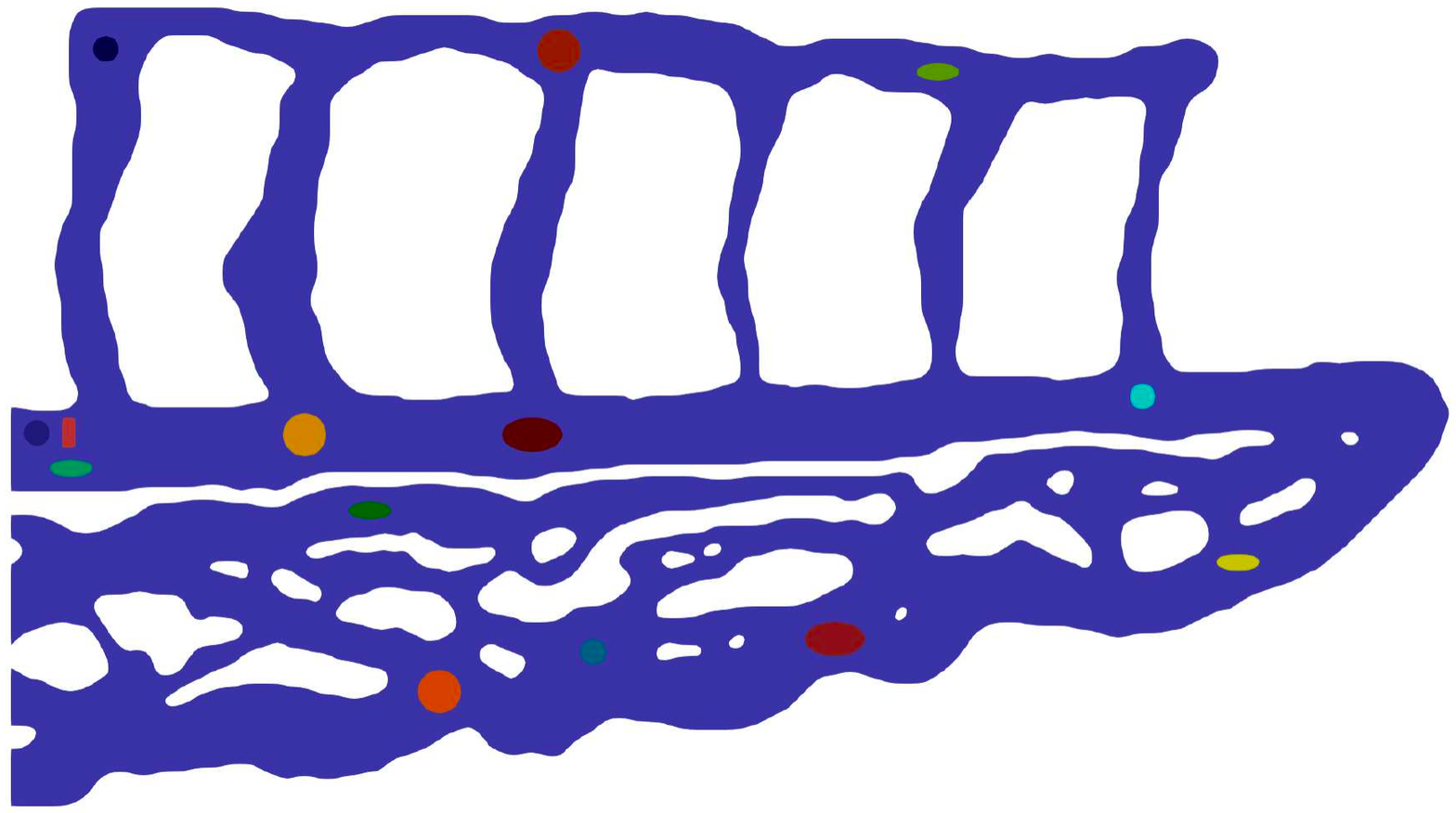}}
    \caption{Creation of two-dimensional bodies inside a complex discretized environment.}
    \label{insertion2D}
\end{figure}

\begin{figure}[h!]
    \begin{center}
		\begin{subfigure}{0.495\textwidth}
			\centering
			\includegraphics[scale=0.18]{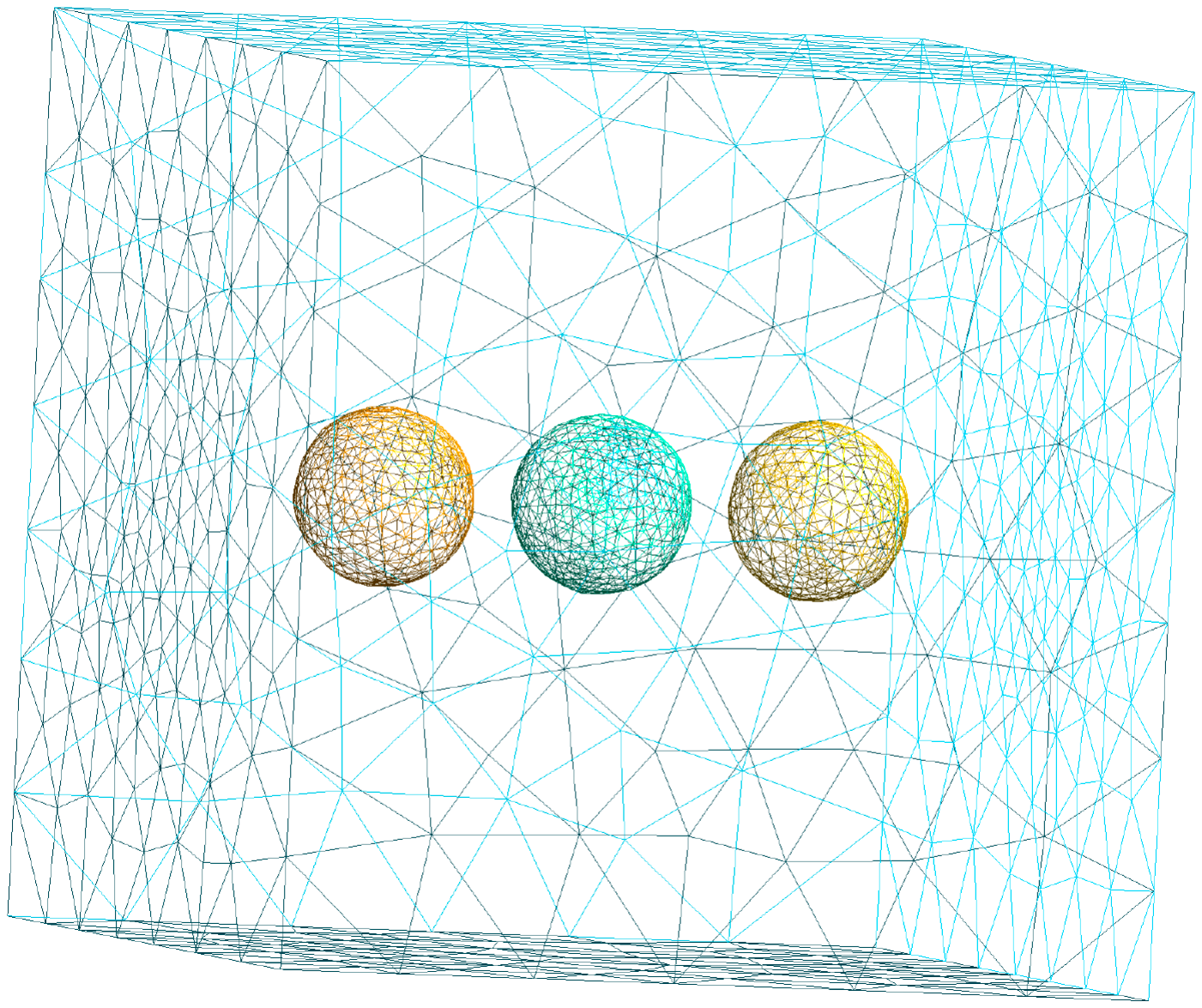}
			\caption{Three-sphere swimmer}
		\end{subfigure}
		\begin{subfigure}{0.495\textwidth}
			\centering
			\includegraphics[scale=0.18]{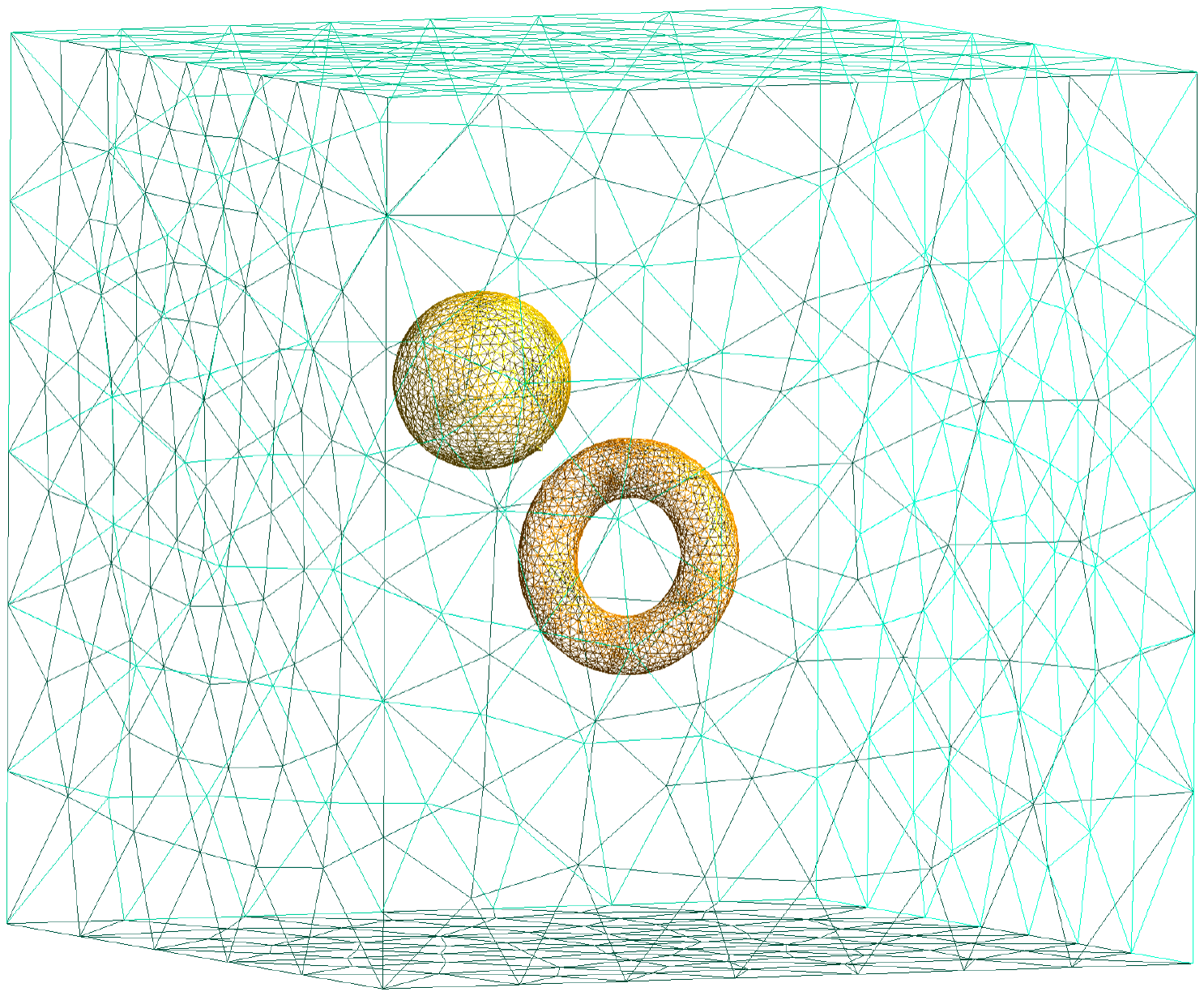}
			\caption{Complex shaped bodies}
		\end{subfigure}
		
		\end{center}

    \caption{Two cases showing the insertion of bodies in three dimensions.}
    \label{insertion3D}
\end{figure}

The capabilities of the tool are illustrated in figures \ref{insertion2D} and
\ref{insertion3D}. Figure \ref{insertion2D} shows a two-dimensional 
environment of complex geometry where a collection of objects have been 
successfully inserted without changes 
to their shapes. In particular all corners of the rectangular body have 
been recovered. Figure \ref{insertion3D} illustrates a simpler 
three-dimensional domain where multiple bodies have been again 
successfully inserted. All three resulting meshes can be used for the 
simulation of fluid-solid interaction with collision treatment, 
as described in the first part of the paper.

\section*{Acknowledgements}

This work of the Interdisciplinary Thematic Institute IRMIA++, as part of the ITI 
2021-2028 program of the University of Strasbourg, CNRS and Inserm, was supported 
by IdEx Unistra (ANR-10-IDEX-0002), and by SFRI-STRAT’US project (ANR-20-SFRI-0012) 
under the framework of the French Investments for the Future Program.
The authors acknowledge the financial support of the French Agence Nationale de 
la Recherche (grant ANR-21-CE45-0013 project NEMO), and 
Cemosis. In addition, the authors thank A. Froehly and L. Cirrottola
for their helpful comments and discussions.

\bibliography{paper_contact}

\end{document}